\documentclass[a4paper,12pt]{amsart}
\DeclareRobustCommand{\SkipTocEntry}[5]{}
\usepackage[english]{babel}
\usepackage[T1]{fontenc}
\usepackage{array}
\usepackage{makecell}
\usepackage[a4paper,margin=3cm]{geometry}
\usepackage[justification=centering]{caption}

\usepackage{amsmath,amsthm,amssymb,amsxtra,calligra,mathrsfs}
\usepackage{graphicx}
\usepackage{tikz}
\usepackage{tikz-cd} 
\usepackage{xcolor}
\usepackage{upgreek}
\usepackage{comment}
\usepackage{graphicx}
\usepackage{mathtools}
\usepackage{extarrows}
\usepackage{faktor}
\usepackage{makecell}
\usepackage{enumitem}
\usepackage{comment}
\usepackage{bookmark}
\usepackage{hyperref}
\usepackage{mathrsfs}
\usepackage{multirow}
\hypersetup{
    colorlinks=true, 
    citecolor=blue,
    linkcolor=magenta,
    pdfauthor={Haoyu Wu},
	bookmarksnumbered=true,
}

\newcommand{\bbF}{\mathbb{F}}

\newcommand{\bG}{\mathbb{G}}

\newcommand{\sheafHom}{\mathscr{H}\text{\kern -3pt {\calligra\large om}}\,}

\newcommand*{\rom}[1]{\expandafter\@slowromancap\romannumeral #1@}

\DeclareMathOperator{\Proj}{Proj}

\DeclareMathOperator{\Ext}{Ext}

\DeclareMathOperator{\Hilb}{\mathrm{Hilb}}

\DeclareMathOperator{\PGL}{\mathrm{PGL}}
\DeclareMathOperator{\GL}{\mathrm{GL}}
\DeclareMathOperator{\SL}{\mathrm{SL}}

\DeclareMathOperator{\Pic}{\mathrm{Pic}}
\DeclareMathOperator{\NS}{\mathrm{NS}}
\DeclareMathOperator{\vol}{\mathrm{vol}}
\DeclareMathOperator{\ord}{\mathrm{ord}}
\DeclareMathOperator{\NE}{\mathrm{NE}}
\DeclareMathOperator{\NL}{\mathrm{NL}}
\DeclareMathOperator{\sh}{\mathrm{Sh}}

\newcommand{\II}{I\!I}
\newcommand{\bA}{{\mathbb A}}
\newcommand{\bZ}{{\mathbb Z}}

\newcommand{\bC}{{\mathbb C}}
\newcommand{\bQ}{{\mathbb Q}}
\newcommand{\bP}{{\mathbb P}}
\newcommand{\bR}{{\mathbb R}}

\newcommand{\fX}{{\mathfrak X}}

\DeclareMathOperator{\wt}{wt}

\DeclareMathOperator{\val}{Val}

\DeclareMathOperator{\pic}{Pic}
\DeclareMathOperator{\aut}{Aut}

\DeclareMathOperator{\sog}{SO}
\DeclareMathOperator{\coeff}{coeff}
\DeclareMathOperator{\mdiag}{diag}

\newcommand{\q}{/\!\!/}
\newcommand{\cO}{{\mathcal O}}
\newcommand{\cX}{{\mathcal X}}

\newcommand{\cP}{{\mathcal P}}
\newcommand{\cv}{{\mathcal V}}
\newcommand{\cY}{{\mathcal Y}}

\newcommand\proj{\text{\rm Proj}}

\newcommand{\cF}{{\mathcal F}}

\newcommand\cD{{\mathcal{D}}}
\newcommand\cS{{\mathcal{S}}}

\newcommand\cE{{\mathcal{E}}}

\newcommand\cL{\mathcal{L}}

\newcommand\cC{{\mathcal{C}}}

\newcommand\rP{\mathrm{P}}

\newtheorem{defn-pro}{Definition-Proposition}
\newtheorem{defn-thm}{Definition-Theorem}
\newtheorem{thm}{Theorem}[section]
\newtheorem{lem}[thm]{Lemma}
\newtheorem{cor}[thm]{Corollary}
\newtheorem{pro}[thm]{Proposition}

\newtheorem{defn}[thm]{Definition}

\newtheorem{question}[thm]{Question}

\newtheorem{conj}[thm]{Conjecture}
\newtheorem{rem}[thm]{Remark}
\newtheorem{case}{Case}

\theoremstyle{remark}

\title{K-moduli of log Del Pezzo pairs}
\author{Long Pan}
\address{Fudan University}
\email{18110180008@fudan.edu.cn}

\author{Fei Si}
\address{Beijing International Center for Mathematical Research, Peking University, No. 5 Yiheyuan Road Haidian District, Beijing, P.R.China 100871}
\email{sifei@bicmr.pku.edu.cn}

\author{Haoyu Wu}
\address{Fudan University}
\email{hywu18@fudan.edu.cn}
\date{}

\begin{document}
\begin{abstract}
We establish the full explicit wall-crossing for K-moduli space $\overline{P}^K_c$ of degree $8$ del Pezzo pairs $(X,cC)$ where  generically $X \cong \bbF_1$ and $C \sim -2K_X$. We also show  K-moduli spaces $\overline{P}^K_c$ coincide with Hassett-Keel-Looijenga(HKL) models $\cF(s)$ of a $18$-dimensional  locally symmetric spaces associated to the lattice $E_8\oplus U^2\oplus E_7\oplus A_1$ under the transform $s(c)= \frac{1-2c}{56c-4}$. Some discussions with relation to  KSBA moduli spaces are also provided.
\end{abstract}

\maketitle
\tableofcontents
\section{Introduction}
\subsection{Background}
This is a continuation of our investigation on the birational geometry of moduli spaces $P_d$ of  del Pezzo pairs of degree $d$. Here a del Pezzo pair $(X,C)$ of degree $d$ consists of a del Pezzo surface $X$ with $(-K_X)^2=d$ and a curve $C\sim -2K_X$. In \cite{PSW1}, we propose to study the Hasset-Keel-Looijenga (HKL) program for  $P_d$  to connect various compactifications of $P_d$. Moreover, the HKL program for $P_d$ should coincide with the wall crossings of K-moduli spaces $\overline{P}^K_{d,c}$.

In this paper, we focus on the K-moduli space $\overline{P}^K_c=\overline{P}^K_{8,c}$  compactifications of del Pezzo pairs of degree $8$ and confirm the above proposal in the case of degree $8$. 
The method should work for other toric del Pezzo pairs. 

Recall for a del Pezzo pair $(\bbF_1, C)$,  one can obtain a K3 pair $(Y,\tau)$  where $Y \rightarrow X$ is a doubld covering branched along curve $C$ and  $\tau: Y\rightarrow Y$ is the natural involution. Then the complement of  Néron-Severi group $\NS(Y)$ in K3 lattice has  orthogonal lattice $\Xi=U^2\oplus E_7 \oplus E_8 \oplus A_1$.   Thus  one can associate a period point in the locally symmetric variety $ \cF:= \widetilde{O}(\Xi) \setminus \cD_\Xi$ for degree $8$ del Pezzo pair $(\bbF_1,C)$ where  $ \cD_\Xi$ is the period domain of $\Xi$. In particular, there is period map  \[ p: P_8 \rightarrow \cF\]
The Torelli theorem for the pairs implies $P_8\subset \cF$ is an open subset and thus $P_8$ has Baily-Borel compactification $P^\ast=\cF^\ast$. 
The lattice  embedding  $ L \hookrightarrow \Xi=U^2\oplus E_7 \oplus E_8 \oplus A_1 $ will induce the natural morphism of locally symmetric variety  \begin{equation}\label{shmor}
    \sh(L):= \widetilde{O}(L) \setminus \cD_L   \longrightarrow \cF
\end{equation} The Noether-Lefschetz locus $\NL(L)$ on $\cF$ is defined as the image of natural morphism (\ref{shmor}). In particular, there are two specific codimension one Noether-Lefschetz locus known as Heegener divisors. One may refer to \cite{BLMM} for general theory of Noether-Lefschetz locus in moduli spaces of K3 surfaces. These two  divisors are defined by the complement lattice  \[ L=(E_7\oplus A_2)_{\II}^\perp \ \ \hbox{and}\ \ L=(E_8 \oplus A_1)_{\II}^\perp\]  in the  Borcherds lattice $\II=E_8^3\oplus U^2$, which we call hyperelliptic divisor $H_h$ and unigonal divisor $H_u$.   
In \cite{PSW1},  we give an arithmetic stratification for the Baily-Borel compactification $P^\ast$ by using sublattice tower in  Borcherds lattice $\II=E_8^3\oplus U^2$. More precisely, looking at  the complement  $(\Xi)^\perp_{\II} \cong E_7\oplus A_1$    and  by adding roots and some modifications,  \cite[Section 5]{PSW1} 
provide three types of the tower for hyperelliptic divisor $H_h$
\begin{enumerate} 
    \item  $A$ type:   $\ E_7\oplus A_1 \subset  E_7\oplus A_2 \subset  E_7\oplus A_3 \subset  \cdots  \subset E_7\oplus A_5$ and its modification tower \[ (E_7\oplus A_5)'\subset (E_7\oplus A_6)'\subset (E_7\oplus A_6)'' \subset (E_7\oplus A_7)'' \subset (E_7\oplus A_7)'''\]
    where the modified lattice $(E_7\oplus A_n)'$ is spanned by $E_7\oplus A_n$ and a root $\delta$, similarly $(E_7\oplus A_n)''$ is spanned by $E_7\oplus A_n$ and two roots \footnote{These roots could be long roots.} $\delta_1, \delta_2$ and so on for  $(E_7\oplus A_n)'''$
    \item  $D$ type: $\ E_7\oplus D_3=E_7\oplus A_3 \subset  E_7\oplus D_4 \subset  \cdots \subset  E_7\oplus D_{6}$ and its modification tower \[ (E_7\oplus D_7)'\subset (E_7\oplus D_8)'\subset  (E_7\oplus D_8)'' \subset  (E_7\oplus D_8)''' \]
    \item  $E$ type:  $\  E_7\oplus E_6 \subset  E_7\oplus E_7 \subset  E_7\oplus E_8$
\end{enumerate}
The complement of the lattice from  the  above tower will produce natural lattice embedding $L\hookrightarrow \Xi $. 
Thus it makes sense to 
denote  $\NL(A_n)$ the Noether-Lefschetz locus associated with the complement of $E_7\oplus A_n$. 
The similar notation $\NL(A_n')$ for the modified lattice in the tower. Then there is stratification for $H_h$ given by 
\begin{equation*}
  \begin{split}
       A \hbox{-type}: & \NL(A_4) \subset \NL(A_3)   \subset \NL(A_2) \subset  \NL(A_1)=\cF^\ast \\
        A\hbox{-modified type}: & \NL(A_7'') \subset \NL(A_7') \subset \NL(A_6')   \subset \NL(A_5') \subset  \NL(A_5)\\
       D \hbox{-type}: &\NL(D_7) \subset  \cdots   \subset \NL(D_4) \subset  \NL(D_3)=\NL(A_3) \\
      D\hbox{-modified type}: &  \NL(D_8') \subset    \NL(D_7') \subset  \NL(D_7)\\
      E \hbox{-type}:&   \NL(E_8) \subset \NL(E_7) \subset \NL(E_6) \subset \cF^\ast
  \end{split}  
\end{equation*}

Similarly,  \cite[Section 5]{PSW1} also introduces towers of latices for unigonal divisor $H_u$ as follows
\begin{equation*}
    \ E_7\oplus A_1 \subset  E_8\oplus A_1 \subset  E_8\oplus A_1^2 \subset E_8  \oplus A_2^2  \subset   E_8 \oplus D_4 \oplus A_2
\end{equation*}
where the rank jump number two is due to a  modification.  
Denote $\NL(U_n)$ the Noether-Lefschetz locus associated with the complement lattice $(E_8 \oplus A_n)^\perp_{\II}$ as before. Then an arithmetic stratification for $H_u$  is obtained as follows
\begin{equation*}
    \hbox{unigonal type}:  \NL(U_4'') \subset \NL(U_3') \subset \NL(U_2') \subset \NL(U_1)=H_u.
\end{equation*}
 Define the  scheme
\begin{equation}
\cF(s):=\Proj(R(\cF^\ast, \Delta(s))), \ \ s\in [0,1] \cap \bQ
\end{equation}
where the $\bQ$-line bundle $\Delta(s)$ is given by \[ \Delta(s)=\lambda+s(H_h+25H_u) .\]
Here $\lambda$ is the Hodge line bundle on locally symmetric variety $\cF^\ast$. $\cF(s)$ is called Hassett-Keel-Looijenga (HKL) model. 
In \cite{PSW1}, according to the arithmetic strategy  under the assumption that $ R(\cF^\ast, \Delta(s))$ is finitely generated, the walls for $\cF(s)$ when  $s\in [0,1] \cap \bQ$ varies is predicted as 
\begin{equation}\label{arwall}
    \{ \ \frac{1}{n}\ |\ n=1,2,3,4, 6,8,10,12,16,25,27,28 ,31\ \}.
\end{equation}

 Let  $\overline{P}^{GIT}$ be the GIT partial compactification space of smooth del pezzo pairs (see Section \ref{section1stwall}  for precise definition) birational to $\cF^\ast$.
The  HKL program  for degree $8$ del Pezzo pair  proposed in \cite{PSW1} is the following conjecture. 
\begin{conj}[HKL for $P_8$] \label{HKL8}
 Notation as above, 
\begin{enumerate}
    \item The section rings $ R(\cF^\ast, \Delta(s))$ are finitely generated for all $s$. In particular, $\cF(s)$ is a projective variety of dimension $18$. 
    \item  $\cF(s)$ will interpolate $P^\ast $ and $\overline{P}^{GIT}$. 
    \item There is isomorphism $\overline{P}_{c}^K \cong \cF(s)$ under the transformation 
    \begin{equation*}
        s=s(c)= \frac{1-2c}{56c-4}.
    \end{equation*}
    In particular, the K-moduli walls coincide with HKL walls given by (\ref{arwall}).
\end{enumerate} 
\end{conj}
The main purpose of the paper is to very the above conjectural picture for resolution of birational period map $\overline{P}^{GIT} \dashrightarrow \cF^\ast$.  In particular, we will provide a resolution of this birational period map with modular meanings. 

\subsection{Main results}
Our first results establish the full wall crossings for the K-moduli spaces $\overline{P}^K_c$. Let $\overline{P}^K_c$ be the good moduli space of K-semistable pairs of degree $8$. Denote $\pi :X \rightarrow \bP^2$ or  $\pi :X \rightarrow \bP(1,1,4)$ the blowup at a point $p$. For any curve  $C$ on $X$ such that  $C \sim -2K_X$, we write $B=\pi(C)$. Let $Y$ be the double cover of $X$ branched along $C$ if $X=\bbF_1$ or the double cover of minimal resolution $\widetilde{X}$ of $X$ branched along the proper transform $\widetilde{C}$ of curve $C$ if $X=Bl_{[1,0,0]}\bP(1,1,4)$.  
\begin{thm}\label{mainthm1} Notation as above, then
\begin{enumerate}
    \item The walls 
    for K-moduli can be divided into  hyperelliptic type walls  $W_h$ and unigonal type walls  $ W_u$ where
    \begin{equation}
    \begin{split}
        W_h=&\{\ \frac{1}{14},\frac{5}{58},\frac{1}{10},\frac{7}{62},\frac{1}{8},\frac{5}{34},\frac{1}{6},\frac{7}{38},\frac{1}{5},\frac{5}{22} ,\frac{2}{7} \ \} \\
        W_u=& \{\  \frac{29}{106},\frac{31}{110},\frac{2}{7}, \frac{35}{118} \ \}
    \end{split}    
    \end{equation}
 Moreover,  the center $Z_w$ of each wall for the K-moduli space $\overline{P}_c^K$ is either a point or a rational curve. The curve  parametrized by $Z_w$ are listed in the Table \ref{Kwall1} and Table \ref{Kwall2}
    \begin{center}
\renewcommand*{\arraystretch}{1.2}
\begin{table}[ht]
    \centering
      \begin{tabular}{ |c  |c |c|c|}
    \hline
     wall &  curve $B$ on $\bP^2$ & weight  & curve singularity at $p$  \\ \hline 
     
     $\frac{1}{14} $  &   $x^4zy=0$  & (1,0,0)  &   $A_1$
     \\ \hline
     
     $\frac{5}{58}$ & $x^4z^2+x^3y^3=0$ & (0,2,3)  & $A_2$
     \\ \hline
     
     $\frac{1}{10}$ & $x^4z^2+x^3zy^2+a\cdot x^2y^4=0,a\in \bC^\ast$ & (0,1,2) &  $A_3$
      \\ \hline
      
     $\frac{7}{62}$ & $x^4z^2+xy^5=0$ & (0,2,5)
     &   $A_4$
     \\ \hline
     
     \multirow{2}{*}{$\frac{1}{8}$} & $x^4z^2+x^2zy^3+a\cdot y^6=0,\ a\in \bC^\ast$  & (0,1,3)  &  $A_5$ tangent to  $L_z$
     \\ \cline{2-4}
     & $x^3f_3(z,y)=0$  & (0,1,1)  &  $D_4$ 
     \\ \hline
     
     \multirow{2}{*}{$\frac{5}{34}$ }
     & {$x^4z^2+xzy^4=0$}  & (0,1,4)
     & $A_7$ \hbox{with a line}  \\ \cline{2-4}
     & $x^3z^2y+x^2y^4=0$ & (0,2,3)
     & $D_5$  
     \\ \hline
     
     \multirow{2}{*}{$\frac{1}{6}$ }
     & {$x^4z^2+zy^5=0$} & (0,1,5)
     &  $A_9$ \hbox{with a line}  \\ \cline{2-4}
     &  $x^3z^2y+x^2zy^3+a\cdot xy^5=0,a\in \bC^\ast$ & (0,1,2)
     &  $D_6$ 
     \\ \hline
     
     \multirow{2}{*}{$\frac{7}{38}$}
     &  $x^3z^2y+y^6=0$ & (0,2,5) &
      $D_7$ tangent to $L_z$ \\ \cline{2-4}
     & $x^3z^3+x^2y^4=0$& (0,3,4)
     &  $E_6$
     \\ \hline
     
     $\frac{1}{5}$ & $x^3z^2y+xzy^4=0$ &(0,1,3)
     &  $D_8$  \hbox{with }  $L_z$ 
     \\ \hline
     
     \multirow{2}{*}{$\frac{5}{22}$ }
     & {$x^3z^2y+zy^5=0$}  & (0,1,4)
     & $D_{9}$  \hbox{with } $L_z$ \\ \cline{2-4}
     & $x^3z^3+x^2zy^3=0$ & (0,2,3)
     & $E_7$ 
     \\ \hline
     
     $\frac{2}{7}$ &    $x^3z^3+xy^{5}=0$   &(0,3,5) &  $E_8$ 
     \\ \hline
     
   \end{tabular}
    \caption{ K-moduli walls from  Gorenstein del Pezzo $\bbF_1=Bl_{[1,0,0]}\bP^2$}
    \label{Kwall1}
\end{table}
\end{center}
where  $L_z$  is the line defined by $\{ z =0\}$ in the Table \ref{Kwall1}.  The weight in the Table \ref{Kwall1} or \ref{Kwall2} means  the weight of $\bG_m$-action on $\bP^2$ or $\bP(1,1,4)$.

\begin{center}
\renewcommand*{\arraystretch}{1.2}
\begin{table}[ht]
    \centering
      \begin{tabular}{ |c |c |c|c|c|c|}
    \hline
     wall &  curve $B$ on $\bP(1,1,4)$ & weight  & $(a,b,m)$ & singularity\\  \hline
      $\frac{29}{106}$   &  $z^3+z^2x^4=0$   &  (1,0,4) & $(0,1,0)$  & $A_1$\\ \hline 
     $\frac{31}{110}$  &  $z^3+zyx^7=0$   &(2,0,7)  &$(1,1,1)$ & $A_1$ with a  tangent line  \\ \hline   
     $\frac{2}{7}$ &     $z^3+y^2x^{10}=0$   &(3,0,10)  & $(2,1,2)$ &  $A_2$ with a  tangent line \\ \hline  
     
$\frac{35}{118}$ &    $z^3+zy^2x^6+y^3x^9=0$ & (1,0,3)  & (1,0,1)  & $D_4$  with a  tangent line \\ \hline 
   \end{tabular}
    \caption{ K-moduli walls from  index $2$ del Pezzo  $Bl_{[1,0,0]}\bP(1,1,4)$}
    \label{Kwall2}
\end{table}
\end{center}

\item  There are birational morphisms 
 \[ \overline{P}^K_{w+\epsilon}\  \xrightarrow{p^+}\  \overline{P}^K_{w}  \  \xleftarrow{p^-}  \  \overline{P}^K_{w-\epsilon}\] 
 for each $w\in W_h\cup W_u$. If $w= \frac{5}{58}$ (rep.  $w=\frac{29}{106}$), $p^+$ is a divisorial contraction with exceptional divisor birational to hyperelliptic divisor $H_h$ (resp. unigonal divisor $H_u$) and $p^-$ is an isomorphism. 
 For the remaining walls $w$, $p^+$ and $p^-$ are flips. Let $E_w^\pm$ be the exceptional locus of $p^\pm$, then $E_w^+$ are described in the Table \ref{tabAE}, \ref{tabDE}, \ref{tabE}, \ref{tabu} and  $E_w^-$ are described in the Table \ref{F1-}, \ref{P-}. 
 In the following tables, $Q$ is an irreducible plane quintic curve. 
\end{enumerate}

\begin{center}
\begin{table}[ht]
    \centering
      \begin{tabular}{ |c | c |c |c|c|}
    \hline
     wall $w$   & Branched curve  $B$ &  $B$ vs. line  $L_x$   &  NL loci  & $\dim \NL$  \\ \hline 
      $\frac{5}{58}$  &        irreducible, $A_2$ at $p$   &  $L_x\cap B=\{3p,p',p'',p'''\}$  & $\NL(A_2)$ & 17 \\ \hline 
     $\frac{1}{10}$    &   irreducible ,  $A_3$  at $p$  &  $L_x\cap B=\{4p,p',p''\}$  & $\NL(A_3)$ & 16 \\ \hline  
     $\frac{7}{62}$   &     irreducible,   $A_4$ at $p$  &  $L_x\cap B=\{5p, p' \}$  & $\NL(A_4)$ & 15 \\ \hline 
    $\frac{1}{8}$    &  irreducible, $A_5$ at $p$ & $B$  tangent  to  $L_x$   & $\NL(A_5')$   & 13 \\ \hline 
    $\frac{5}{34}$  &   $L_x+Q$,   $A_7$ at  $p$  &  $L_x\cap Q=\{4p,p'\}$  & $\NL(A_6'')$  &  11  \\ \hline 
    $\frac{1}{6}$ &    $L_x+Q$, $A_9$  at  $p$  &  $L_x\cap Q=\{5p\}$  &  $\NL(A_7''')$    &  9 \\ \hline 
   \end{tabular}
    \caption{Geometric descriptions of curves in exceptional locus $E^+_{w}$ of $A$-type wall-crossing on $\bbF_1$}
    \label{tabAE}
\end{table}
\end{center}

\begin{center}
\begin{table}[ht]
    \centering
      \begin{tabular}{ |c | c |c |c|c|}
    \hline
     wall $w$  & Branched curve  $B$ &  $B$ vs. line  $L_x$   &  NL loci & $\dim \NL$ \\ \hline 
     $\frac{1}{8}$ &  irreducible,  $D_4$ at $p$ & $L_x\cap B=\{3p,p',p'',p'''\}$    & $\NL(D_4)$ & $15$ \\ \hline 
    $\frac{5}{34}$ &    irreducible ,  $D_5$ at $p$&  $L_x\cap B=\{4p,p',p''\}$   & $\NL(D_5)$& $14$ \\ \hline 
    $\frac{1}{6}$ &      irreducible,  $D_6$   at $p$&  $L_x\cap B=\{5p,p'\}$  &  $\NL(D_6)$ & $13$ \\ \hline 
    $\frac{7}{38}$ &       irreducible, $D_{7}$  at $p$ &  $L_x\cap B=\{6p\}$  &  $\NL(D_{7}')$ & $11$ \\ \hline
     $\frac{1}{5}$ &   $L_x+Q$, $D_{8}$  at $p$  &  $L_x\cap Q=\{4p,p'\}$  &  $\NL(D_{8}')$ & $10$ \\ \hline
     $\frac{5}{22}$ &    $L_x+Q$,   $D_{10}$ at  $p$  &  $L_x\cap Q=\{5p\}$  &  $\NL(D_{9}')$  &  $9$ \\ \hline
   \end{tabular}
    \caption{Geometric description of curves in exceptional locus $E^+_{w}$ of $D$-type wall-crossing on $\bbF_1$}
    \label{tabDE}
\end{table}
\end{center}

\begin{center}
\begin{table}[ht]
    \centering
      \begin{tabular}{ |c | c |c |c|c|}
    \hline
     wall $w$    & Branched curve  $B$ &  $B$ vs. line  $L_x$   &  NL loci & $\dim \NL$  \\ \hline 
    $\frac{7}{38}$ &      irreducible,   $E_6$ at $p$  &  $L_x\cap B=\{4p,p',p''\}$  &  $\NL(E_6)$   &  $13$\\ \hline 
     $\frac{5}{22}$  &   irreducible,   $E_7$ at $p$  &  $L_x\cap B=\{5p,p'\}$ &  $\NL(E_{7})$ &  $12$ \\ \hline
     $\frac{2}{7}$ & irreducible,   $E_8$ at $p$  & $L_x\cap B=\{5p,p'\}$  &  $\NL(E_{8})$  &  $11$ \\ \hline
   \end{tabular}
    \caption{Geometric description of curves in exceptional locus $E^+_{w}$ of $E$-type wall-crossing on $\bbF_1$}
    \label{tabE}
\end{table}
\end{center}

\begin{center}
\begin{table}[ht]
    \centering
      \begin{tabular}{ |c |c |c |c|}
    \hline
     wall $w$ &    curve  $B$ on  $\bP(1,1,4)$   & NL loci &   $\dim \NL$ \\ \hline 
      $\frac{29}{106}$  &    irreducible,  $A_1$ at $p$   &  $\NL(U_1)$ & $17$  \\ \hline 
      
     $\frac{31}{110}$ & \makecell[c]{  irreducible, $A_1$ at $p$ and  \\ tangent to the line $L_y$  and \\
     passing through another fixed point}   & $\NL(U_2')$ &  $15$  \\ \hline  
     $\frac{2}{7}$   &  \makecell[c]{  irreducible,  $A_2$ at $p$ \\ tangent to the line $L_y$ and  \\   passing through another fixed point} &     $\NL(U_3’)$  &  $14$ \\ \hline 
     $\frac{35}{118}$  &  \makecell[c]{ irreducible,  $D_4$ at $p$ \\  tangent to the line $L_y$ and \\  passing through another fixed point}  & $\NL(U_4'')$    &   $12$ \\ \hline 
   \end{tabular}
    \caption{Geometric description of curves in  exceptional locus $E^+_w$ of K-moduli wall-crossing on $X=Bl_p\bP(1,1,4)$}
    \label{tabu}
\end{table}
\end{center}
\end{thm}

\begin{rem}
    For $d=9$, the wall crossing results are established  in \cite{ADL19}.
    For $d=8$ and the surface is deformation equivalent to $\bP^1 \times \bP^1$, the explicit wall crossing  is known due to \cite{ADL2021}.  In \cite{ADL2021}, the authors proved  K-moduli spaces of $(X,cC)$ pairs is isomorphic to a global VGIT for all $c\in (0,\frac{1}{2})$ and the walls follow from the computation of Laza-O'Grady \cite{LaO18} directly. But in our case it seemss there is no global VGIT construction. 
\end{rem}
\begin{rem}
   After the double cover construction, locally  if  the curves    $B$ on $\bP^2$ or  $\bP(1,1,4)$ has  $ADE$ singularity, then so  does $Y$. In particular, $Y$ is a K3 surface possibly with ADE singularities. In the table above, the Noether-Lefchetz locus $\NL$  parameterizes such K3 surface $Y$ obtained by the curve $B$. 
\end{rem}

Our second result is to establish the following isomorphism of K-moduli space and HKL models.
\begin{thm}[= Theorem \ref{thm-identification}]\label{mainthm2}
    There is natural isomorphism $\overline{P}_c^K \cong \cF(s)$induced by the period map   under the transformation  \[ \ s=s(c)= \frac{1-2c}{56c-4}\] where $\frac{1}{14}<c<\frac{1}{2}$. In particular, $P_c^K$ will interpolates the GIT space $\overline{P}^{GIT}$ and Baily-Borel compactification $\cF^\ast$.  
\end{thm}
Combining the results of theorem \ref{mainthm1} and theorem \ref{mainthm2}, we also deduce the following directly
\begin{thm}
     Conjecture \ref{HKL8} holds.
\end{thm}

As a application of our results, we can also determine some walls for K-moduli space of Fano 3-fold pairs via cone construction (\cite[Theorem5.2]{ADL22}).
\begin{cor}[= Corollary \ref{3foldpair}]
The stability threshold of Fano 3-fold pairs $(X,cS)$  from cone construction of del Pezzo pairs   $(\bbF_1,w_nC)$
 are given by the rational numbers
\begin{equation*}
   \{\ c=\frac{11+n}{27+n} \ |\  n=1,2,\cdots, 5\ \} \cup \{\ c=\frac{3+n}{11+n}\ |\ n=6, 7,8,9,11\  \} 
\end{equation*}
Similarly,  the stability threshold of Fano 3-fold pairs $(X,cS)$ from $(Bl_p \bP(1,1,4),w_nC)$ are given by
\begin{equation*}
     \{\ c=\frac{36+m}{52+m}\ |\ m=1,3,4,7\  \}.
\end{equation*}
In particular, the above rational numbers are walls for K-moduli space of Fano $3$-fold pairs. 
\end{cor}

\addtocontents{toc}{\SkipTocEntry}
\subsection*{Relations to other works}
 \subsubsection*{(1)} In \cite{ABBDILW23}, the authors construct a compact moduli spaces $P_d^{CY}$ parameterizing  S-equivalence classes of log Calabi-Yau pairs  $(\bP^2, \frac{3}{d}C_d)$ and their degenerations. Their work builds a birdge connecting K-moduli of plane curves in \cite{ADL19}  and KSBA moduli of plane curves by Hacking \cite{Hacking04}. 
It would be  interesting to extend their work to our setting.  Note that for the KSBA moduli space parametrizing  del pezzo surface pairs $(X,cC)$ with $c=\frac{1}{2}+\epsilon$,  \cite{AEH}  has related  it to the toroidal compactification of locally symmetric varieties.  

 \subsubsection*{(2)} Our results also provide the birational contractability of Heenger divisors.  It would be useful to detect walls for other  HKL model on a locally symmetric variety,  as long as the locally symmetric variety is associated to  a lattice $\Lambda$ of signature  $(2,n)$ and admitting an embedding $E_8\oplus U^2\oplus E_7\oplus A_1 \hookrightarrow \Lambda$ (up to a lattice saturation).  For example, in the forthcoming work   \cite{FLLST},  the authors study the HKL on locally symmetric variety given by lattice $E_8\oplus U^2\oplus E_7\oplus A_2$, where they use the birational contractability of Heegner divisors obtained as birational transform of excetional locus in our K-moduli $\overline{P}_c^K$ wall-crossing  to give prediction the walls for HKL of modui space of quasi-polarised K3 of genus $4$.

\addtocontents{toc}{\SkipTocEntry}
\subsection*{Organization of the paper}
In section  \ref{sect2},  we collect basic results from algebraic K-stability theory, including  the equivariant K-stability criterion, K-moduli theory and the existentce of wall-crossings in K-moduli. In Section \ref{sect3},  we discuss the K-degeneration of del Pezzo pairs based on local to global volume comparison and $T$-singularities theory. In Section \ref{sect4} we use the method of equivariant K-stability criterion from complexity one del Pezzo pairs to find all walls for K-moduli space $\overline{P}^K_c$. Based on the computation in Section \ref{sect4}, we finish the explicit wall crossing in Section \ref{explicitwc}.  In Section \ref{sect6}, we realize the K-moduli space $\overline{P}^K_c$ as the HKL model under the transformation of the parameter. 
In the final section \ref{sect7}, we discuss some possible  relations on the K-moduli space $\overline{P}^K_c$ to  the moduli of log CY pairs and  KSBA moduli theory for the pairs.

\addtocontents{toc}{\SkipTocEntry}
\subsection*{Acknowledgement}

We would like to thank  Zhiyuan Li and Yuchen Liu for many helpful discussions. 
The second author would also like to thank Chen Jiang for the helpful discussion on quotient singularities and Junyan Zhao for the comments on the drafts. 
Part of this work was written when the second author visited SCMS and he would like to thank their hospitality.
The second author was partially supported by LMNS (the Laboratory of Mathematics for Nonlinear Science, Fudan University).
This project is supported by the NKRD Program of China (No.2020YFA0713200), NSFC Innovative Research Groups (Grant No.12121001) and General Program (No.12171090).

\section{Preliminaries on K-stability and K-moduli } \label{sect2}
\subsection{Notation of K-stability} All varieties are over $\bC$. 
A pair $(X,D)$ is consisting of a projective normal variety $X$ with an effective $\bQ$-divisor $D$ such that $K_X+D$ is $\bQ$-Cartier. By taking a log resolution $f: Y \rightarrow (X,D)$, we have $$K_Y+D_Y=f^\ast(K_X+D).$$  We call a pair $(X,D)$ is log canonical (lc) if  $D_Y$ has coefficients $\leq 1$ and Kawamata log terminal (klt) if $D_Y$ has coefficients $< 1$. 

\begin{defn}
A pair $(X,D)$ is called log Fano if  $(X,D)$ is klt and  the $\bQ$-Cartier divisor $-K_X-D$ is ample.
\end{defn}

Let $E$ be a prime divisor on a birational model $\pi: Y \rightarrow X$ of $X$. The divisor $E$ defines a valuation $v_e$ on the function field $\bC(X)$, which is known as divisorial valuation.  Given such a divisorial valuation  $v_E$, define its $A$-function 
\[   A_{(X,D)}(E):=1+ \coeff _E(K_Y-\pi^\ast(K_X+D)) \]
and $S$-function
\begin{equation*}
    S_{(X,D)}(E):=\frac{1}{(-K_X-D)^n} \int_0^\tau \vol(-\pi^\ast(K_X+D)-tE) dt
\end{equation*} 
where $$\tau:= \sup \{  t \, |\,  -\pi^\ast(K_X+D)-tE\ \hbox{is pseudo-effective}  \ \}$$  is the pseudo-effective threshold. 

\begin{defn-thm}[Fujita-Li] \label{FL} Notation as above,
\begin{enumerate}
    \item The pair $(X,D)$ is K-semistable if and only if \[\beta _{(X,D)}(E):= A_{(X,D)}(E)- S_{(X,D)}(E) \geq 0\] for any prime divisor  $E$ over $X$ 
    \item The pair $(X,D)$ is uniform K-stable if and only if \[ \delta (X,D):= \mathop{\inf } \limits_{E}  \frac{A_{(X,D)}(E)}{S_{(X,D)}(E)}  > 1\]
    where $E$ runs over all prime divisors  over $X$.
\end{enumerate}
\end{defn-thm}

\subsection{Normalised Volume }

Let $\val_{X,x}$ be the valuation space centered at the point $x\in X$. Motivated by the problem of {K\"ahler}-Einstein metric, Chi Li in \cite{MR3715806} introduced the notation of normalised volume, which plays a very important role in controlling singularities of K-semistable Fano varieties.

\begin{defn}[Chi Li] Let $x\in (X,D)$ be a klt singularity. Define the normalised volume function $\widehat{\vol}_{(X,D),x}: \val_{X,x} \rightarrow \bR_{>0} \cup \{\infty \}$ by 
\begin{equation*}
    \widehat{\vol}_{(X,D),x}(v)= \begin{cases}
 A_{(X,D)}(v)^n \cdot \vol(v)& \text{if }\ A_{(X,D)}(v) < \infty ,\\
\infty  & \text{if }\   A_{(X,D)}(v) = \infty . \end{cases}
\end{equation*}
The normalised volume of klt singularity $x\in (X,D)$ is defined as 
\begin{equation*}
    \widehat{\vol}(X,D;x):= \mathop{ \inf } \limits_{v\in  \val_{X,x}}   \widehat{\vol}_{(X,D),x}(v)
\end{equation*}
\end{defn}

\begin{thm}(Local-to-Global volume comparison \cite{MR3872852}) \label{volcom}
 Let $(X,D)$ be a K-semistable log Fano pair, then
\begin{equation*}
    (-K_X-D)^n \leq (1+\frac{1}{n})^n \cdot \widehat{\vol}(X,D;x)
\end{equation*}
for any closed point $x\in X$.
\end{thm}

\begin{thm}(Finite degree formula \cite{xuzhuang21} ) \label{finite}
Let $f: (Y,D_Y) \rightarrow (X,D_X)$ be a quasi-\'etale morphism, then 
\begin{equation*}
    \widehat{\vol}(Y,D_Y;y)=\deg (f) \cdot \widehat{\vol}(X,D_X;x)
\end{equation*}
\end{thm}

\subsection{$\bG_m$-equivariant  K-stability} In general, determining whether a given log Fano pair is K-semistable is a challenging problem in algebraic K-stability theory (see \cite{xu21}). 
But for log Fano pairs of complexity $1$, due to the work of Zhuang and Ilten-Süß, there is an effective method to detect K-stability. 
It is a main tool for us to find walls of K-moduli spaces. Now let's briefly recall the theory. 
Let $T$ be a maximal torus of $\aut(X,D)$. 
By a deep theorem of \cite{JHHLX}, reductivity of  $\aut(X, D)$ is a necessary condition for $(X, D)$ to be K-polystable, so we may assume this and thus $T$ is the unique (up to conjugation) maximal torus of  $\aut(X, D)$.  
Denote $\bC(X)^T$ the $T$-invariant rational function on $X$. Assume $T$ acts on $(X,D)$ effectively. 
\begin{defn}
We say a $n$-dimensional pair $(X,D)$ is of complexity  $m$  where $m=n-\dim T$.
\end{defn}
Clearly, $X$ is of complexity $0$ if and only if $X$ is a toric  variety. 
In this paper, we are interested in  the surface of complexity $1$. 
Such surface $X$ admits a natural rational map $f: X \dashrightarrow \bP^1$ with $\bC(X)^T=\bC(\bP^1)$.  
\begin{defn}
A prime divisor $E$ on  $X$  is called $T$-vertical if the maximal $T$-orbit in $E$ has dimension $\dim T$. Otherwise, it is called  $T$-horizontal.
\end{defn}
Geometrically, $T$-vertical divisors on $X$ can be viewed as the  fibers of the natural rational map $f: X \dashrightarrow \bP^1$ .
The following criterion for log Fano pairs of complexity $1$ is useful. 
\begin{thm}(\cite{IS17}, see also \cite[Theorem 3.2]{liu23})\label{thmcomplexity1}
Let $(X,D)$ be a log Fano pair of dimension $2$ with an effective $\bG_m$-action $\lambda$.  Then $(X,D)$ is K-polystable if and only if the following conditions hold:
\begin{enumerate}
    \item $\beta_{(X,D)}(F)> 0$ for all vertical $\lambda$-invariant prime divisors $F$ on $X$;
    \item  $\beta_{(X,D)}(F)= 0$ for all horizontal $\lambda$-invariant prime divisors $F$ on $X$;
    \item $\beta_{(X,D)}(v)=0$ for the valuation $v$ induced by the 1-PS $\lambda$.
\end{enumerate}
\end{thm}

\subsection{K-moduli stack and its good moduli space}
For $c\in (0,\frac{1}{2})\cap \bQ$,  we consider the moduli stack of pairs $\cP^K_c\colon  Sch(\bC) \rightarrow Set$ by 
\begin{equation*}
\begin{split}
    \cP^K_c(B)=\big\{ & (\cX,\cD) \rightarrow B \ |  \ \hbox{smoothable K-semistable \ family } \\ &  \hbox{with}\  \cD_b \sim_\bQ -2c K_{\cX_b}\ \hbox{for\ any}\  b\in B \ \big\}.
\end{split}
\end{equation*}
The framework of \cite{ADL19} \cite{LXZ} establishes the existence of good moduli space $\overline{P}_c^K$ in the sense of J. Alper \cite{MR3237451}. 
\begin{thm} The moduli stack
$\cP^K_c$ is a separated Artin stack admitting a  good moduli space $\overline{P}_c^K$
$$\cP^K_c \rightarrow \overline{P}_c^K$$
which is a projective normal scheme. Moreover, $\overline{P}_c^K(\bC)$ parametrizes $S$-equivalence of $c$-K-polystable log Fano pairs.
\end{thm}
When we vary the coefficient $c\in (0,\frac{1}{2})\cap \bQ$, there are wall-crossing phenomenon that are established in \cite{ADL19}(see also \cite{zhou} for more explanations). More precisely, they prove
\begin{thm}\label{kwc}
There are finitely many rational numbers (i.e., walls ) $0<w_1<\cdots <w_m <\frac{1}{2}$ such that 
\[ \overline{P}_c^K\cong \overline{P}_{c'}^K \ \  \hbox{for\ any}\ w_i<c,c'<w_{i+1}\ \hbox{and \ any }\ 1 \leq i \leq m-1 .\] Denote $\overline{P}_{(w_i,w_{i+1})}^K:=\overline{P}_c^K$ for some  $c\in (w_i,w_{i+1})$, then at each wall $w_i$ there is a flip (or divisorial contraction)
\begin{equation*}
    \begin{tikzcd}
\overline{P}_{(w_{i-1},w_{i})}^K \arrow[dr,"p^-"] & & \overline{P}_{(w_i,w_{i+1})}^K \arrow[dl,"p^+"']\\
  & \overline{P}_{w_i}^K &
\end{tikzcd}
\end{equation*}
which fits into a local variation of geometric invariant theory (VGIT) in the sense of \cite[Section 2.2]{AFS17}.
\end{thm}
At a wall $w$, let $U_w \subset \overline{P}_w^K$ be the open subset where $p^\pm$ is isomorphic. Denote $Z_w:=\subset \overline{P}_w^K-U_w$ the complement and $E_w^\pm:=(p^\pm)^{-1}(Z_w)$. Locally $p\pm|_{E_w\pm}:\ E_w^\pm \rightarrow Z_w$ is a fibration with fiber a weighted projective space by the general theory of VGIT (see \cite{DH}) and local VGIT interpretation of K-moduli space. 

\begin{defn-pro}(see \cite[Section 2.4]{ADL19}) The CM-line bundle $\lambda_{c}(f)$ on $Z_{c,m}^{red}$ is defined by   
\begin{equation*}
    c_1(\lambda_{c}(f))=-f_\ast((K_f+c\cD)^3)
\end{equation*}
where the pushforward is in the sense of cycle.
\end{defn-pro}

\begin{thm}(\cite[Theorem1.1]{XZ20} or \cite[Theorem 3.36]{ADL19})\label{cmample} The $\lambda_{c}(f)$ can descend to the good moduli space  $\overline{P}^K_c$ and it is ample. 
\end{thm}

\section{Moduli of del Pezzo pairs and K-semistable  degenerations }\label{sect3}

\subsection{ K-polystable degeneration  of log Fano pairs}
Let $(\cS,c \cD) \rightarrow (B,o)$ be a  flat and $\bQ$-Goresntein family of log Fano pairs over a smooth pointed curve $(B,o)$. Assume generic fiber $(\cS_\eta,c \cD_\eta)$ is a K-semistable smooth pair and set $(X,c C)=(\cS_o,c D_o)$, it is a very important and also difficult problem in K-moduli theory to understand the geometry of the degeneration $(X,cC)$. In general, we have
\begin{thm}(Odaka \cite{MR3010808} for absolute version)
If $(\cS_0,c \cD_o)$ is K-semitable, then $(\cS_o,c \cD_o)$ has klt singularity at worst. In particular, $\cS_o$ is normal.
\end{thm}
Recall  the ADE surface singularities  are isolated singularities defined by the following local equations
\begin{equation*}
    \begin{split}
        & A_n: x^2+y^2+z^{n+1}=0, \\
      & D _n: x^2+y^2z+z^{n-1}=0, n\ge 4 \\
       & E_6:  x^2+y^3+z^4=0, \\
       &  E_7:  x^2+y^3+yz^3=0, \\
       &  E_8:  x^2+y^3+z^5=0. \\
    \end{split}
\end{equation*}

By the classification of klt surface singularities (see \cite{MR3057950}),   $(\cS_0,c D_0)$ has quotient singularities. We are only interested in such quotient singularities that are smoothable and these singularities  are so called $T$-singularity in the literature (see \cite{Hacking04} \cite{HP10}).

\begin{defn}($\bQ$-Gorenstein smoothing)
Given  a log Fano pair $(X,c D)$, we call a proper flat morphism $(\cX,\cD) \xrightarrow{f} B$ is a $\bQ$-Gorenstein smoothing of $(X,D)$ if 
\begin{enumerate}
   \item $-K_{\cX/\cD}$ is $\bQ$ Cartier-divisor and $f$-ample
    \item $\cD\sim -rK_{\cX/\cD}$ is an effective $\bQ$ Cartier-divisor on $\cX$ and $supp(\cD)$ does not contain any fiber. Here the index $0<r<$.
    \item $f$  and $f|_\cD$ is smooth over $B-\{o\}$ and $(\cX_o,\cD_o) \cong (X,D)$.
\end{enumerate}
\end{defn}

In the dimension $2$ without boundaries, $T$-singularity can be classified.
\begin{pro} \label{class}
A $T$-singularity is one of the following 
\begin{enumerate}
    \item $A_n$, $D_n$ or $E_6,\ E_7,\ E_8$;
    \item  Cyclic Quotient singularity of type $\frac{1}{ln^2}(1,lna-1)$ for some $l,n,a\in \bZ_{\geq 0}$ and $\gcd(a,n)=1$, i.e.,  $(0\in \bA^2_{x,y}/\mu_{ln^2})$ given by $$ \xi \cdot (x,y)=(\xi \cdot x,\xi^{l na-1} \cdot y )$$
    where $\xi$ is $ln^2$-th primitive root.
\end{enumerate}
\end{pro}
\begin{proof}
See \cite[ Proposition 3.10]{HP10}.
\end{proof}
We collect basic properties of $T$-singularity as follows. 
\begin{pro} \label{basic} Let $X$ be a projective surface with $T$-singularity at worst.
\begin{enumerate}
    \item  the following Noether type formula holds 
    \begin{equation*} 
        K_X^2+e_{top}(X)+\mathop{\sum} \limits_{p\in X_{sing}} \mu_p=12\chi(\cO_X)
    \end{equation*}
    where $\mu_p=b_2(M_p)$ is the second Betti number of   the Milnor fiber $M_p$ of its smoothing family. In particular, \begin{equation*}
\mu_p= \begin{cases}
 r & \text{if }\ p\  is\  A_r,D_r\  or\  E_r,\\
l-1& \text{if }\  p \ is\ of\  type\  \frac{1}{ln^2}(1,lna-1).
\end{cases}
\end{equation*}
Moreover, if $X$ is rational, then 
\begin{equation}\label{Noeth}
 K_X^2+\rho(X)+\mathop{\sum} \limits_{p\in X_{sing}} \mu_p =10   
\end{equation}
\item if in addition $-K_X$ is big and nef, then $$\dim H^0(X,-mK_X)=\frac{m(m+1)}{2}K_X^2+1.$$
\end{enumerate}
\end{pro}

Based on the above characterization of $T$-singularity, we will give a Cartier index  estimate for the K-semistable degeneration of del Pezzo surface pairs of degree $d=8$. 
\begin{thm} \label{Cartier index}
Let $(X,C)$ be the a central fiber of a family $(\cX,\cC) \rightarrow B$ over a curve $B$ where the general fiber $(\cX_b,\cC_b)$ is a smooth pair  of degree $d=8$ , then  the Cartier index of canonical bundle satisfies $ind(K_X,x) \leq 3$.
\end{thm}
\begin{proof}
By the deformation invariance of $K_X^2$ and  the formula (\ref{Noeth}), the T-singularities appeared on $X$ must be of type  \[\frac{1}{n^2}(1,an-1)\ \hbox{or}\   \frac{1}{2n^2}(1,2an-1).\]
By the local to global volume comparison (\ref{volcom}) and finite degree formula (\ref{finite}), 
\begin{equation} \label{vol}
    \frac{4d}{9}(1-2c)^2 \leq  \widehat{\vol}(X,c D;x)=\frac{1} {n^2}  \widehat{\vol}(\widetilde{X},c \widetilde{D};\widetilde{x}) \leq  \frac{1} {n^2}(2-c \cdot \ord_{\widetilde{x}}(\widetilde{D}))^2
\end{equation}
where $(\widetilde{X},\widetilde{D}) \rightarrow (X,D)$ is a smooth covering. By the Skoda's inequality\[2-c \cdot \ord_{\widetilde{x}}(\widetilde{D})>0.\] Thus, if $\ord_{\widetilde{x}}(\widetilde{D}) \geq 4$, then the inequality (\ref{vol}) will show
\begin{equation*}
    \frac{2n \sqrt{d}}{3} \leq 2-(\ord_{\widetilde{x}}(\widetilde{D})-4)\frac{c}{1-2c} \le 2,
\end{equation*}
which implies $n \leq \frac{3}{\sqrt{d}}$. Thus, $n=1$.

Now we may assume $1 \leq \ord_{\widetilde{x}}(\widetilde{D}) =i+j \leq 3$ where  $x^iy^j$ is the monomial with minimal degree of defining equation for $\widetilde{D}$ under coordinate $(x,y)$ for $\widetilde{X}$. As $2K_{\widetilde{X}}+\widetilde{D}$ descends to a Cartier divisor $2K_X+D \sim 0$, then 
\begin{equation} \label{indef}
    2an \equiv i+(na-1)j\ \mod n^2
\end{equation}
In particular, $ i \equiv j \mod\ n$ holds. 

If $i+j=3$, then by $n | (i-j)$, we know $n=1$ or $n=3$. For the latter case, if $(i,j)=(0,3)$, the indefinite equation (\ref{indef})  is just \[an \equiv 3 \ \mod\ n^2,\]  thus $(a,n)=(1,3)$ gives a solution.  

If $i+j=2$, $n | i-j$ shows $n=1$ or $2$.  This finishes proof.
\end{proof}

Combining the classification results in \cite[Table 6]{Nakayama} and \cite[Table 10]{fujitaYasutake17}, we have 
\begin{cor}\label{classification of degeneration}  Let $(X,C)$ be the K-semistable degeneration of log del Pezzo pairs of degree $d=8$,   
\begin{enumerate}
    \item if the Cartier index of $K_X$ is $2$, then  $X$ is isomorphic to the blowup along a smooth point of $\bP(1,1,4)$; 
    \item if the Cartier index of $K_X$ is $3$, then  $X$ is the surface in  Subsection \ref{index3surface}.
\end{enumerate}
\end{cor}
\begin{proof}
 (1) is  the direct consequence of \cite[Table 6]{Nakayama} and \cite[Proposition 7.1]{Nakayama}. (2) follows from classification results in \cite{fujitaYasutake17}.
\end{proof}
\subsection{Geometry of log del Pezzo of index $3$}\label{index3surface}

Let $\bbF_n$ be the Hirzebruch surface and $\sigma,f$ be the section and fiber. Thus, we have intersection numbers
\[ \sigma^2=-n,\ f^2=0,f.\sigma=1\]
According to Nakayama's strategy (see  \cite{Nakayama}) of  classification of log del Pezzo surface  $X$ of index $a>1$,  one can first take a minimal resolution $\pi: M \rightarrow X$ to obtain a nonsingular rational surface $M$ and thus there is a birational morphism \[\mu: M \rightarrow S , \  S \ \hbox{is either}\  \bbF_n \ \hbox{or} \ \bP^2. \]
Then one inductively runs MMP to decompose  $\mu$ into a series of birational  morphisms contracting $(-1)$-curves. It turns out that  one can recover the surface $X$ from $M$ and the dual graph   of exceptional divisors of  $\pi$. In the case $a=3$ and degree $d=8$,  Fujita-Yasutake's results (see \cite[Table 10]{fujitaYasutake17}) show that $S=\bbF_5$  and $\mu$ is just blowup of distinct two points on the fiber $f$ but not in the section $\sigma$ with exceptional divisor $E_1$ and $E_2$. Then \[ \pi:\ M \longrightarrow X \]
  is the morphism contracting  two  divisors $F_1, F_2$ 
  where $F_1$ and $F_2$ are proper transforms of section $\sigma$ and fiber $f$  so  that
   \[ F_1^2=-5,\ F_1.F_2=1,\ F_2^2=-2.\]
 Thus, $F_1$ and $ F_2$ are contracted to a singularity of type $\frac{1}{9}(2,1)$ and
 \begin{equation*}
     -3K_\pi=-3(K_M-\pi^\ast K_X)=2\cdot F_1+F_2.
 \end{equation*}
 We can write the intersection matrix on $M$ as follows:
 \begin{equation}\label{M} \left(
    \begin{array}{c|c c c c }
    &F_1 & F_2 & E_1 & E_2  \\
    \hline
    F_1 & -5 & 1 & 0 & 0\\
    
    F_2 & 1 & -2 & 1 & 1\\ 
    
    E_1 & 0 & 1 & -1 & 0\\
    
    E_2 & 0 & 1 & 0 & -1 \\
    
    \end{array}
  \right) \end{equation}  

\vspace{0.3cm}
\section{Computation of K-moduli walls}\label{sect4}
\subsection{$S$-function on $\bbF_1$ and  $Bl_p\bP(1,1,4)$}  
To apply equivariant K-stability criterion  in Theorem \ref{thmcomplexity1}, we need to compute $S$-function associated to exceptional divisor of weighted blowup centered at the invariant points under the maximal torus action. First, observe that the boundary $C\in  |-2K_X|$ will imply  
\begin{equation}\label{scaling}
\begin{split}
     S_{(X,cC)}(F) &=\frac{1}{(-K_X-cC)^2}\int^{\infty}_{0}vol(-\mu^*(K_X+cC)-tF)dt\\
    &=\frac{1-2c}{(-K_X)^2}\int^{\infty}_{0}vol(-\mu^*K_X-tF)dt,
\end{split}
\end{equation}
where the second identity is due to the change of variable. Thus, the computations of $S(F)$-function for divisor $F$ is reduced to the computation  of volume of  divisor $L_t:=-\mu^*K_X-tF$. In the computation of volume, Zariski decomposition on surface is a very useful tool and let's recall the following 
\begin{pro}\label{zariski}
If $X$ is a normal projective surface and $D$ is an pesudo-effective $\bQ$- divisor on $X$, then there is a unique decomposition 
\[ D=P+N \]
where $P, N$ are two $\bQ$- effective divisors such that $P.N_i=0$ for each irreducibel component of $N$, $P$ is nef and the intersection matrix of irreducible components of $N$ is negative or $N=0$. In particular, \[ \vol(D)=P^2.\]
\end{pro}

Second, we recall the basics of weighted blowup in dimension 2 following \cite[Chapter 3.2]{prok}. Let $X=(\bC^2/\bZ_r(a,b),0)$ be the local germ of 2-dimensional quotient singularity of order $r$ where $a,b,r\in \bZ_{>0}$ with $\gcd(a,b)=1$. A weighted blowup $\pi: Y \rightarrow X$ of weight $(a,b)$ is  a projective birational morphism which is isomorphic over $X-\{0\}$ and $E:=\pi^{-1}(0)$.
 By \cite[Lemma 3.2.1]{prok}, we have
 \begin{equation} \label{wbu}
 \begin{split}
     K_Y=&\pi^\ast K_X+(-1+ \frac{a+b}{r})E \\
     \pi^\ast C=&\widetilde{C}+m_C(a,b)E 
 \end{split}
 \end{equation}
where $\widetilde{C}$ is the proper transform of $C$ and 
\[ m_C(a,b):= \min \{ a\cdot i+b \cdot j\ |\  x^iy^j\  \hbox{is a monomial of local equation  defining } C   \ \}.\]
\subsubsection{$S$-functions of $\bG_m$-invariant prime divisors on $\bbF_1$}\label{4.1.1}
Let $H_x,H_y,H_z$ be the proper transform of divisors $\{x=0\}$,$\{y=0\}$ and $\{z=0\}$  under the blowup $\mu:\bbF_1\rightarrow  \bP^2 $ of $\bP^2$ at point $[0,0,1]$. Then 
\begin{equation*}
    H_x \sim H_y \sim H_z-E \sim \mu^\ast \cO_{\bP^2}(1) -E
\end{equation*}
where $E$ is exceptional divisor. Moreover,  
\begin{equation*}
    H_x^2=0,\ H_x.E=1.
\end{equation*}
By a theorem of Miyaoka (see also \cite[Lemma2.1]{fulger11}), the nef cone $Nef(\bbF_1)$ of $\bbF_1$ is generated by $H_x$ and $E+H_x$, or equivalently, the Mori cone  $\NE(\bbF_1)$  of $\bbF_1$ is generated by $H_x$ and $E$. Actually, if $C$ is an effective curve on $\bbF_1$, one can write $C \sim a \sigma+b f $ for $a,b\in \bZ$ and $C$ is not zero. Note that 
 \begin{equation*}
    p_\ast (\cO_{\bbF_1}(a \sigma+b f))=\left\{
    \begin{aligned}
    0 & , & \  a <0, \\
   \cO_{\bP^1}(b) \otimes sym^a(\cO\oplus \cO(-1)) & , & a \ge 0.
\end{aligned}
\right.
\end{equation*}
implies  $a \ge 0$ and  $b \ge 0$. Conversely, 
\[ H^0(\cO_{\bbF_1}(a \sigma+b f)) \ne 0 \ \ \hbox{if}\ \  a \ge 0,\ b \ge 0 \]

\begin{pro} \label{sfunction1}
Notation as above, then
\begin{equation*}
\begin{split}
      S_{(\bbF_1,cC)}(H_z)=&\frac{5}{6}(1-2c),\ \ \
      S_{(\bbF_1,cC)}(H_x)=S_{(\bbF_1,cC)}(H_y)=\frac{13}{12}(1-2c)\\
      S_{(\bbF_1,cC)}(E)=&\frac{7}{6}(1-2c).
\end{split}
\end{equation*}
\end{pro}
\begin{proof}
Clearly, we have 
\[ -K_X =3H_x+2E . \]
To compute the volume function , it is sufficient  to find Zariski decomposition for $L_t:=-K_X-t F$ where
$F$ is  either $E$, $H_x$, $H_y$ or $H_z$.
\begin{enumerate}
    \item If $F=E$, then \[L_t= 3H_x+(2-t)E=(2-t)(E-H_x)+(5-t)H_x\] is nef for $t \le 2$ and 
    \[ L_t^2=(4+t)\cdot (2-t)\] 
    Thus, the psudo-effective threshold  is $\tau=2$ and thus 
    \begin{equation*}
        \begin{split}
         \int_0^{2}  \vol(L_t) dt=&\int_0^{2}  L_t^2 dt  =\frac{28}{3} 
        \end{split}
    \end{equation*}
    and the formula (\ref{scaling}) shows $S_{(\bbF_1,cC)}(E)=\frac{7}{6}(1-2c)$.
    \item If $F=H_x$, then $L_t= (3(1-2c)-t)H_x+2(1-2c)E$ is nef for $ 0 \le  t\le (1-2c) $ and  thus its volume is
    \[L_t^2= -4(1-2c)t+8(1-2c)^2 .\] 
    For $ (1-2c) \le t$, the Zariski decomposition  $L_t=P_t+N_t$ where
    \[P_t=(3(1-2c)-t) (H_x+E)\]
    Therefore, we get 
    \[P_t^2=(3(1-2c)-t)^2\]
    In particular, the pesudo-effective threshold is $3(1-2c)$ and thus 
    \begin{equation*}
        \begin{split}
           \int_0^{\tau}  \vol(L_t) dt=& \int_0^{(1-2c)}L_t^2 dt+ \int_{1-2c}^{3(1-2c)} P_t^2 dt 
           = \frac{26}{3} (1-2c)^3
        \end{split}
    \end{equation*}
    \item If $F=H_z$, then $L_t= (3(1-2c)-t)H_x+(2(1-2c)-t)E$ is nef for $ 0 \le  t$ and \[ L_t^2=(4(1-2c)-t)\cdot (2(1-2c)-t)\]
    Thus, the pesudo-effective threshold is $2(1-2c)$ and
    \begin{equation*}
        \begin{split}
           \int_0^{\tau}  \vol(L_t) dt=& \int_0^{3(1-2c)}L_t^2 dt
           = \frac{20}{3} (1-2c)^3
        \end{split}
    \end{equation*}
\end{enumerate}
Then we finish the proof by applying the formula (\ref{scaling}).
\end{proof}
Fix the  blowup $\pi:\bbF_1\cong Bl_p\bP^2 \rightarrow\bP^2$ where $p=[1,0,0]$. For complexity one pair $(\bbF_1,cC)$, we can choose a suitable coordinate $[x,y,z]$ on $\bP^2$ such that the maximal torus action $\lambda:\bG_m\rightarrow \aut(\bbF_1,C)$ can be regarded as the lifting of the action $\lambda(t)[x,y,z]=[t^{\lambda_1}x,t^{\lambda_2}y,t^{\lambda_3}z]$.
To compute $S$-function of the divisorial valuation determined by $\lambda$, we need to figure out the relation between the weight $(\lambda_1,\lambda_2,\lambda_3)$ and the induced valuation $\ord_F$.
\begin{case}\label{Case1}
\begin{equation*}
    \min\{\lambda_1,\lambda_2,\lambda_3\}=\lambda_2 (\text{resp}.~ \lambda_3)
\end{equation*}
\end{case}
In this case  $\ord_F$ is given by the exceptional divisor of the weighted blowup along the point $[0,1,0]$ (\text{resp}.~ $[0,0,1]$) with weight
    \[\wt(x)=\lambda_1-\lambda_2,\wt(z)=\lambda_3-\lambda_2\] 
    \[(\text{resp}.~\wt(x)=\lambda_1-\lambda_3,\wt(y)=\lambda_2-\lambda_3)\]
    on the local parameter $(x,z)$( (\text{resp}.~ $(x,y)$). To simplify the symbol, we write 
    \[
    a:=\lambda_1-\lambda_2, b:=\lambda_3-\lambda_2
    \]
    \[
    (\text{resp}.~ a:=\lambda_1-\lambda_3, b:=\lambda_2-\lambda_3)
    \]

\begin{case}\label{Case2}
\begin{equation*}
   \min\{\lambda_1,\lambda_2,\lambda_3\}=\lambda_1
\end{equation*}
\end{case}  
  In this case the description about $\ord_F$ is slightly complicated. Its center located on the exceptional divisor $E$ over $\bP^2$, but the local parameter not only comes from $\bP^2$. Note that $\bbF_1\subset \bP^2_{[x,y,z]}\times\bP^1_{[u,v]}$ is  defined by equation $\{yv-zu=0\}$ and $\bG_m$-action on $\bbF_1$  is restricted from a $\bG_m$-action $\lambda'$ on $\bP^2\times\bP^1$ given by \[\lambda'(t)([x,y,z],[u,v]):=([t^{\lambda_1}x,t^{\lambda_2}y,t^{\lambda_3}z],[t^{\lambda_2}u,t^{\lambda_3}v]).\] Using the coordinate $[u,v]$, we can accurately describe the induced valuation $\ord_F$. 
    \begin{itemize}
        \item If $\lambda_2>\lambda_3$,  then $F$ is the exceptional divisor of  weighted blowup $Y\rightarrow\bbF_1$ where the center is $([1,0,0],[1,0])$ and  the local parameter is $(z,u)$ with weight \[a:=\wt(z)=\lambda_3-\lambda_1, 
\ b:= \wt(u)=\lambda_2-\lambda_3.\]

        \item If $\lambda_3>\lambda_2$, then $F$ is the exceptional divisor of  weighted blowup $Y\rightarrow\bbF_1$ where the center is $([1,0,0],[0,1])$ and the local parameter is $(y,v)$ with weight
        \[a:=\wt(y)=\lambda_2-\lambda_1,\  b:=\wt(v)=\lambda_3-\lambda_2.\]
    \end{itemize}

We denote all the weighted blowups mentioned above as $\mu\colon Y\rightarrow\bbF_1$. 

\begin{lem}\label{Mori Cone}
If $\mu\colon Y\rightarrow\bbF_1$ is in the Case \ref{Case2} mentioned above, then
\[
\overline{\NE}(Y)=\bR_{\ge0}[F]+\bR_{\ge0}[\overline{E}]+\bR_{\ge0}[\overline{L}]
\]
where $\overline{L}$ is the strict transform of a line $L$ on $\bP^2$ passing through all the blowup centers, and $\overline{E}$ is the strict transform of exceptional divisor $E$ on $\bbF_1$.

\end{lem}

\begin{proof}

It is enough to prove the lemma under condition $\lambda_2>\lambda_3$, the same argument applies to the remaining case. It is easy to check $\NS(Y)$ is generated by $F,\overline{E}$ and $\overline{L}$. Moreover,  their intersection matrix on $Y$ is following
\begin{equation*} \left(
    \begin{array}{c|c c c }
    & F &  \overline{E} & \overline{L} \\
    \hline 
    F & -\frac{1}{ab} & \frac{1}{b}  &  \frac{1}{a} \\
    
    \overline{E} &  \frac{1}{b}  & -1-\frac{a}{b}  & 0\\
   
    \overline{L}  &  \frac{1}{a}  & 0 & -\frac{b}{a} \\
 \end{array}
  \right)  \end{equation*}     
Note that $F,\overline{E}$ and $\overline{L}$ all have negative self intersections, so they are extremal curves by \cite[Section 11.2]{prok}. On the other hand, for any irreducible curve $D\subset Y$, if $D\neq F$ and $\overline{E}$, then we can write the equation of its image on $\bP^2$ by 
\[
x^{d-e}f_e(y,z)+\cdots +xf_{d-1}(y,z)+f_d(y,z)=0.
\]
for some integers $d$ and $e$. Thus the class of $D$ can be written as 
\[
d\overline{L}+(d-e)\overline{E}+(db+(d-e)a-m)F,
\]
where $m=\ord_F(\mu^*(\mu(D))-D)$ is computed as in formula (\ref{wbu}). In the local coordinate $(z,u)$, the equation of $\mu(D)$ is
\[
f_e(z)+\cdots +u^{d-e-1}f_{d-1}(z)+u^{d-e}f_d(z)=0
.\]
$D$ is irreducible and so is $\mu(D)$. Thus there must exists $u^i$ term, which implies that $m\le ib\le db$. So any curve on $Y$ can be written as the positive linear combination of $\{F,\overline{E},\overline{L}\}$, we deduce that $\overline{\NE}(Y)=\bR_{\ge0}[F]+\bR_{\ge0}[\overline{E}]+\bR_{\ge0}[\overline{L}]$.

\end{proof}

\begin{pro} \label{S1} If the weight of $\lambda$ is in the Case \ref{Case1}, we have  
    \begin{equation*}
        S_{(\bbF_1,cC)}(\ord_F)=
        \begin{cases}
          (a+b-\frac{b^2}{12a})(1-2c), &0<b<a\\
          \frac{13a+10b}{12}(1-2c), & b\ge a
        \end{cases}
    \end{equation*}
    
\end{pro}

\begin{pro}\label{S2}
If the weight of $\lambda$ is in the Case \ref{Case2},  we have 
    \begin{equation*}
        S_{(\bbF_1,cC)}(\ord_F)=\frac{14a+13b}{12} (1-2c)
    \end{equation*}
\end{pro}

\begin{proof}
By (\ref{scaling}), we need to compute the volume of $L_t=-\mu^*K_{\bbF_1}-tF$. To do so, it is sufficient to find the positive part of $L_t$ defined in Proposition \ref{zariski}.
Direct computation shows
\begin{equation}\label{intersection number}
    L_t.F=t,\ \ L_t.\overline{E}=1-\frac{t}{b},\ \ L_t.\overline{L}=2-\frac{t}{a}.
\end{equation} According to Lemma \ref{Mori Cone},  $L_t$ is nef when $0\le t\le\min\{2a,b\}$. On the other hand, observe $\sqrt{8ab}>\min\{2a,b\}$ and  \[L_t^2=8-\frac{t^2}{ab}\ge0 \] for $0\le t\le \min\{2a,b\}$.  So we get
 \[ P(L_t)=L_t, \ \ 0\le t\le \min\{2a,b\}  \]
 and we need to compute the positive part $P(L_t)$ for $t\ge\min\{2a,b\}$. We divide the remaining computations into two cases.
 \begin{itemize}
     \item If $\min\{2a,b\}=b$, from (\ref{intersection number}) $ L_t$ fails to be nef for $t>b$ which is due to $L_t.\overline{E}<0$ for $t>b$. So we may assume the negative part for $t>b$ is $s \overline{E}$ for some $s>0$ and thus the positive part $P(L_t)=L_t-s\overline{E}$.  By Proposition \ref{zariski}, we know that \[P(L_t).\overline{E}=(L_t-s\overline{E}).\overline{E}=0,\] which implies $s=\frac{L_t\overline{E}}{(\overline{E})^2}=\frac{t-b}{a+b}$. Thus, the intersection numbers of $P(L_t)$ with $F,\overline{E}$ and $\overline{L}$ are given by
\[ P(L_t).F=\frac{t+a}{a(a+b)},\ \ \ P(L_t).\overline{E}=0,\ \ \ P(L_t).\overline{L}=\frac{2a-t}{a} \]
     which implies that $P_2(L_t)$ is nef $b\le t\le2a$. Then we have
     \[(P(L_t))^2=L_t^2-sL_t\cdot \overline{E} =9-\frac{(t+a)^2}{a(a+b)}> 0,\ \  \  \ \hbox{for}\  b\le t \le 2a. \]
So we need to continue finding the positive part of $P(L_t)$ for $t>2a$.  The trick to find the positive part for $L_t,\  b\le t \le 2a$ will be applied again. We assume the negative part for $P(L_t)$  is 
$s'\cdot \overline{L}$ for some $s'>0$ due to the fact  $P(L_t). \overline{L}<0$ for $t>2a$. Then from \[P(L_t).\overline{L}=(P(L_t)-s'\overline{L}).\overline{L}=0,\] we get $s'=\frac{t-2a}{b}$ and \[ P(L_t)=(3b+2a-t)(\frac{1}{b}\mu^*\pi^*L-\frac{a}{(a+b)b}\mu^*E-\frac{a}{a+b}F).\] Thus, we have the following intersection number 
\[(P(L_t))^2=(L_t-s\cdot \overline{E}-s' \cdot \overline{L})^2= \frac{(3b+2a-t)^2}{(a+b)b}\]
and $P(L_t).F,\ P(L_t).\overline{E},\ P_3(L_t).\overline{L}$ are non-negative for $2a\le t\le2a+3b$.  Therefore,  we conclude that  the pesudo-effective threshold of $L_t$ is \[\tau= 3b+2a.\] In a summary,  the Zariski decomposition of $L_t$ in this case is given by
      \begin{equation*}
    P(L_t)=
    \begin{cases}
      L_t , & \ \ 0\le t\leq b, \\
    L_t-\frac{t-b}{a+b}\overline{E}  , & \ \ b\le t\le2a,\\
   L_t-\frac{t-b}{a+b}\overline{E}-\frac{t-2a}{b}\overline{L}, &\ \  2a\le t\le2a+3b
    \end{cases}
    \end{equation*}
     \item if $\min\{2a,b\}=2a$, the computation is parallel to the case $\min\{2a,b\}=b$ and we leave it to the interested reader. In the case, we have Zariski decomposition
     \begin{equation*}
    P(L_t)=
    \begin{cases}
      L_t , & \  0\le t\leq 2a, \\
    L_t-\frac{t-2a}{b}\overline{L}  , &\  2a\le t\le b,\\
   L_t-\frac{t-2a}{b}\overline{L}-\frac{t-b}{a+b}\overline{E}, &\  b\le t\le2a+3b
    \end{cases}
\end{equation*}
 \end{itemize}
 A direct computation shows that in both cases we have
 \[ \int_{0}^{\infty} \vol(L_t) dt=\int_{0}^{\infty}(P_t)^2 dt=\frac{28a+26b}{3}.  \] So we conclude that \[S_{(\bbF_1,cC)}(F)=\frac{1-2c}{8}\cdot \frac{28a+26b}{3}=\frac{14a+13b}{12}(1-2c).\]

\end{proof}

\subsubsection{S-functions on $Bl_p\bP(1,1,4)$ }\label{4.1.2}

Now we consider S-functions on $Bl_p\bP(1,1,4)$ for $\bG_m$-invariant prime divisors. Let $\pi: X=Bl_p\bP(1,1,4) \rightarrow \bP(1,1,4)$ be the blowup at  the  smooth point $p=[1,0,0]$. As before, denote $H_x,H_y,H_z$ be the proper transform of divisors $\{x=0\}$,$\{y=0\}$ and $\{z=0\}$. Then 
\[ H_z\sim 4H_y+3E,\ \ H_y \sim H_x-E, \ H_x\sim \pi^\ast \cO(1). \]
Their intersection numbers are 
\[H_y^2=-\frac{3}{4},\ H_y.E=1,\ E^2=-1 .\]
Let $C \sim a H_y+b E$ be an effective curve on $Bl_p\bP(1,1,4)$, thus 
\[ H^0(X, \cO_X(C))=H^0(\bP(1,1,4),\cO(a) \otimes \mu_\ast \cO_X((b-a)E)) \ne 0 \] 
provided either $a >0$ and   or $a=0$ and $b\ge 1$. 
\begin{pro}
On  the surface $X=Bl_p\bP(1,1,4)$, we have 
\begin{equation*}
\begin{split}
      S_{(X,cC)}(E)=&\frac{83}{48}(1-2c), \ \ \
      S_{(X,cC)}(H_x)=\frac{41}{24}(1-2c), \\
      S_{(X,cC)}(H_y)=&\frac{53}{24}(1-2c), \ \ \
      S_{(X,cC)}(H_z)=\frac{25}{48}(1-2c)\\
      \end{split}
\end{equation*}
\end{pro}
\begin{proof}
First we note that a non-zero nef $\bQ  $-divisor $L$ on $X$ is of the form 
\begin{equation}\label{claim}
     L \sim _\bQ H_y+ t\cdot  E\ \ \  \hbox{for}\ \ \  \frac{3}{4} \le t \le 1  .
\end{equation}
Then following the procedure of computations for $\bbF_1$, We need to compute the volume function of $L_t:=-K_X-cC-t F$ for
$F$ is  $E$, $H_x$, $H_y$ or $H_z$.  Note that 
\[ -K_X-cC \sim (1-2c)(6H_y+5E).\]
\begin{enumerate}
    \item If $F=E$, then  by the (\ref{claim}) $L_t= 6(1-2c)H_y+(5(1-2c)-t)E$ is nef for $ 0 \le  t\le \frac{(1-2c)}{2}$ and  so its volume is
    \[L_t^2= -t^2-2(1-2c)t+8(1-2c)^2 .\] 
    For $\frac{(1-2c)}{2} \le t$, we get the Zariski decomposition $L_t=P_t+N_t$ of $L_t$ where
    \[P_t= (5(1-2c)-t)(\frac{4}{3}H_y+E).\]
    Therefore, 
    \[ \vol(L_t)=P_t^2=\frac{(5(1-2c)-t)^2}{3}.\]
    In particular, the pesudo-effective threshold is $5(1-2c)$. Thus, we get 
    \begin{equation*}
        \begin{split}
           \int_0^{\tau}  \vol(L_t) dt =&\int_0^{\frac{(1-2c)}{2}}L_t^2 dt+ \int_{\frac{(1-2c)}{2}}^{5(1-2c)} P_t^2 dt
           =\frac{83}{6} (1-2c)^3.
        \end{split}
    \end{equation*}
    \item If $F=H_y$, then $L_t= (6(1-2c)-t)H_y+5(1-2c)E$ is nef for $ 0 \le  t\le (1-2c)$ and 
     \[L_t^2= -\frac{1}{4}(3t-8(1-2c))(t+4(1-2c)).\] 
     For $(1-2c)\ge t$, we get the Zariski decomposition $L_t=P_t+N_t$ of $L_t$
     with 
    \[P_t= (6(1-2c)-t)(H_y+E)\]
     So the pesudo-effective threshold is $6(1-2c)$. Therefore, we have
     \begin{equation*}
        \begin{split}
           \int_0^{\tau}  \vol(L_t) dt 
           =\frac{125+87}{12} (1-2c)^3.
        \end{split}
    \end{equation*}
    \item If $F=H_x$, then $L_t= 6(1-2c)-t)H_y+(5(1-2c)-t)E$ is nef for $ 0 \le  t\le 2(1-2c)$ and 
    \[L_t^2= \frac{1}{4}t^2-3(1-2c))t+8(1-2c)^2.\]
    For $2(1-2c)>t $, the positive part of Zariski decomposition for $L_t$ 
    is \[P_t= \frac{5(1-2c)-t}{3}(4H_y+3E)\] and thus the volume of $L_t$ is 
    \[ P_t^2= \frac{(5(1-2c)-t)^2}{3}.\]
    Therefore, the pesudo-effective threshold is $5(1-2c)$ and 
    \begin{equation*}
        \begin{split}
           \int_0^{\tau}  \vol(L_t) dt =\frac{41}{3} (1-2c)^3
        \end{split}
    \end{equation*}
    \item If $F=H_z$, then $L_t= 6(1-2c)-4t)H_y+(5(1-2c)-3t)E$ is nef for $ 0 \le  t\le (1-2c)$ and 
    \[ L_t^2=3t^2-10(1-2c)t+8(1-2c)^2 .\]
    For $t>(1-2c)$, the positive part of $L_t$ is given by
    \[ P_t=(6(1-2c)-4t)\cdot (H_y+E) \]
    and thus 
    \[ P_t^2=(3(1-2c)-2t)^2 .\]
    So  the pesudo-effective threshold $\tau$ is $\frac{3}{2} (1-2c)$. This gives
      \begin{equation*}
        \begin{split}
           \int_0^{\tau}  \vol(L_t) dt =& \int_0^{1-2c}  L_t^2 dt+  \int_{1-2c}^{\frac{3}{2} (1-2c)} P_t^2 dt
           =\frac{25}{6}(1-2c)^3.
        \end{split}
    \end{equation*}
\end{enumerate}
\end{proof}

Let $\bG_m$ be a 1-PS acting on $\bP(1,1,4)$ with weight $(\lambda_1,\lambda_2,\lambda_3)$. Similar to the case in $\bbF_1$, we can classify the induced divisorial valuation $\ord_F$ by the weight $(\lambda_1,\lambda_2,\lambda_3)$.
In the affine open subset $U=\{x=1\}$,  $\pi^{-1}(U) \subset Bl_{[1,0,0]} \bP(1,1,4)$ is isomorphic to 
\begin{equation*}
    \{ (y,z) \times [u,v] \ |\ yv=zu  \}.
\end{equation*}

\begin{itemize}
    \item[{Case $1'$}\label{unigonal case 1'}] : 
    $\lambda_2-\lambda_1<\lambda_3-4\lambda_1<0$.  Then 
    $F$ is the exceptional divisor of the weighted blow up $\mu:Y\rightarrow X$  where the center is $([1,0,0],[1,0])$ and  the local coordinate is $(z,u)$ with weight \[a:=\mathrm{wt}(z)=-(\lambda_3-4\lambda_1),\ \ b:=\mathrm{wt}(u)=-\lambda_2+\lambda_1+(\lambda_3-4\lambda_1)\]

    \item[{Case $2'$}\label{unigonal case 2'}]: 
    $\lambda_3-4\lambda_1>\lambda_2-\lambda_1>0$. Then
    $F$ is the exceptional divisor of the weighted blowup $\mu:Y\rightarrow X$ where the center is $([1,0,0],[0,1])$ and the local coordinate is $(y,v)$ with weight 
    \[a:=\mathrm{wt}(y)=\lambda_2-\lambda_1,\ \ b:=\mathrm{wt}(v)=\lambda_3-4\lambda_1-(\lambda_2-\lambda_1).\]

    \item[Case $3'$\label{unigonal case 3'}]: 
    $\lambda_1-\lambda_2>0, \lambda_3-4\lambda_2>0$. Then $F$ is the exceptional divisor of the weighted blowup $\mu:Y\rightarrow X$ where the center on $\bP(1,1,4)$ is $[0,1,0]$ and the local coordinate is $(x,z)$ with weight
    \[a:=\mathrm{wt}(x)=\lambda_1-\lambda_2,\ \ b:=\mathrm{wt}(z)=\lambda_3-4\lambda_2.\]
\end{itemize}

The remaining cases can be reduced to the above cases. For example, if $\lambda_1-\lambda_2>0, \lambda_3-4\lambda_2<0$, then it can be reduced to Case1' or Case2' by replace $(\lambda_1,\lambda_2,\lambda_3)$ by $(-\lambda_1,-\lambda_2,-\lambda_3)$.

\begin{lem}\label{Curve Cone}
Let $\mu\colon Y\rightarrow X$ be the weighted blowup mentioned above. We denote the strict transform of $E$, $H_x$, $H_y$, $H_z$ on $Y$ by $\overline{E}$, $\overline{H}_x$, $\overline{H}_y$, $\overline{H}_z$ respectively. Then 
\begin{enumerate}

    \item  in the \hyperref[unigonal case 1']{Case $1'$},  
    \[  \overline{\NE}(Y)=\bR_{\ge0}[F]+\bR_{\ge0}[\overline{E}]+\bR_{\ge0}[\overline{H}_y]\]
    
    \item in the \hyperref[unigonal case 2']{Case $2'$}, 
    \begin{equation*}
    \overline{\NE}(Y)=
        \begin{cases}
          \bR_{\ge0}[F]+\bR_{\ge0}[\overline{E}]+\bR_{\ge0}[\overline{H}_y] & \hbox{if}\ 0<b\le3a \\
          \bR_{\ge0}[F]+\bR_{\ge0}[\overline{E}]+\bR_{\ge0}[\overline{H}_y]+\bR_{\ge0}[\overline{H}_z]  & \ \hbox{otherwise} 
        \end{cases}
    \end{equation*}
    \item in the \hyperref[unigonal case 3']{Case $3'$}, 
    \begin{equation*}
       \overline{\NE}(Y)= 
       \begin{cases}
          \bR_{\ge0}[F]+\bR_{\ge0}[\overline{E}]+\bR_{\ge0}[\overline{H}_x] & 0<b\le3a \\
          \bR_{\ge0}[F]+\bR_{\ge0}[\overline{E}]+\bR_{\ge0}[\overline{H}_x]+\bR_{\ge0}[\overline{H}_z] & a<b\le 4a \\
          \bR_{\ge0}[F]+\bR_{\ge0}[\overline{E}]+\bR_{\ge0}[\overline{H}_z] & b\ge4a
        \end{cases}
    \end{equation*}
   \end{enumerate}
\end{lem}

\begin{proof}
It is clear that $\rho(Y)=3$. Now we determine the Mori cone in the three cases as discussed above.
\begin{itemize}
    \item[{\hyperref[unigonal case 1']{Case $1'$}}]:
    In this case the basis of $\Pic(Y)_{\bR}$ can be chosen by $\{F, \overline{E}, \overline{H}_y\}$. The intersection matrix on $Y$ is the following:
    
   \begin{equation*}
   \left(
       \begin{array}{c|c c c}
    & F & \overline{E} & \overline{H}_y\\
    \hline
    F & -\frac{1}{ab} &  \frac{1}{b} &  \frac{1}{a} \\
    
    \overline{E} &  \frac{1}{b} & -1-\frac{a}{b} & 0\\
    
    \overline{H}_y &  \frac{1}{a}  & 0 & -\frac{3}{4}-\frac{b}{a} \\
    
    \end{array}
    \right)
   \end{equation*}
Note that $F,\overline{E}$ and $\overline{H}_y$ all have negative self-intersections, so they are extremal curves by \cite[Section 11.2]{prok}. On the other hand, for any irreducible curve $D\subset Y$, if $D\neq F$ and $\overline{E}$, then we can write the equation of its image on $\bP(1,1,4)$ by 
    \[
    z^df_l(x,y)+z^{d-1}f_{l+4}(x,y)+...+f_{4d+l}(x,y)=0
    \]
    for some integer $d$ and $l\in\{0,1,2,3\}$. Thus the class of $D$ can be written as 
    \[
    (4d+l)\overline{H}_y+(4d+l-e)\overline{E}+((4d+l)b-m+(4d+l-e)a)F,
    \]
    where $e=\ord_E(\mu(D))$ and $m=\ord_F(\mu^*(\mu(D))-D)$ is computed as in formula (\ref{wbu}). Note that $e\le d+l$. In the local coordinate $(z,u)$, the equation of $\mu(D)$ is
    \[
    z^{-e}(z^{d}f_l(zu)+\cdots +zf_{4d+l-4}(zu)+f_{4d+l}(zu))=0
    .\]
    $D$ is irreducible and so is $\mu(D)$. Thus there must exists $u^i$ term, which implies that $m\le ib\le (4d+l)b$. So any curve on $Y$ can be written as the positive linear combination of $\{F,\overline{E},\overline{H}_y\}$, we deduce that $\overline{\NE}(Y)=\bR_{\ge0}[F]+\bR_{\ge0}[\overline{E}]+\bR_{\ge0}[\overline{H}_y]$.

 \item[{\hyperref[unigonal case 2']{Case $2'$}}] : 
    In this case the basis of $\Pic(Y)_{\bR}$ can be chosen by $\{F, \overline{E}, \overline{H}_y\}$ and  their intersection matrix  is the following
    \begin{equation*}
    \left(
        \begin{array}{ c |c c c  }
    & F& \overline{E} & \overline{H}_y\\
    \hline
    F & -\frac{1}{ab} & \frac{1}{b} & \frac{1}{a}\\
    
    \overline{E} & \frac{1}{b} & -1-\frac{a}{b} & 0\\

    \overline{H}_y &  \frac{1}{a} & 0 & -\frac{3}{4} \\
    \end{array}
    \right)
    \end{equation*}    
Observe that $\{F, \overline{E}, \overline{L}_y\}$ both have negative self-intersection, so they must be extremal. However, $(\overline{H}_z)^2=3-\frac{b}{a}$, its sign depends on $(a,b)$. Using the same arguments as above, we see that
\[
    \overline{\NE}(Y)= 
       \begin{cases}
         \bR_{\ge0}[F]+\bR_{\ge0}[\overline{E}]+\bR_{\ge0}[\overline{H}_y] & \ \hbox{if}\  0<b\le3a \\
        \bR_{\ge0}[F]+\bR_{\ge0}[\overline{E}]+\bR_{\ge0}[\overline{H}_y]+\bR_{\ge0}[\overline{H}_z] &  \ \hbox{if}\  b>3a \\
        \end{cases}
\]

\item[{\hyperref[unigonal case 3']{Case $3'$}}] :    
    In this case the basis of $\Pic(Y)_{\bR}$ can be chosen by $\{F, \overline{E}, \overline{H}_x\}$. The intersection matrix on $Y$ is the following
    \begin{equation*}
    \left(
        \begin{array}{c|c c c }
    & F & \overline{E} & \overline{H}_x \\
    \hline
     F  & -\frac{1}{ab}  & 0 &  \frac{1}{a} \\
    
    \overline{E} &  0  & -1 & 0\\
    
    \overline{H}_x &  \frac{1}{a}  & 0 & \frac{1}{4}-\frac{a}{b} \\
    
    \end{array}
    \right)
    \end{equation*}
We can see that $F$ and $\overline{E}$ are extremal. Note that
    \[(\overline{H}_x)^2=\frac{1}{4}-\frac{a}{b},\ \ (\overline{H}_z)^2=3-\frac{b}{a} .\] 
    Hence
    \begin{itemize}
        \item if $0<b\le3a$, then $\overline{H}_x$ is extremal but $\overline{H}_z$ is not extremal.
        
        \item if $3a<b\le4a$, then $\overline{H}_x$ and $\overline{H}_z$ are both extremal.
        
        \item if $b\ge4a$, then $\overline{H}_z$ is extremal but $\overline{H}_x$ is not extremal.
    \end{itemize}
    
    Using the same arguments again we conclude that

    \begin{equation*}
       \overline{\NE}(Y)= 
       \begin{cases}
          \bR_{\ge0}[F]+\bR_{\ge0}[\overline{E}]+\bR_{\ge0}[\overline{H}_x] & 0<b\le3a \\
         \bR_{\ge0}[F]+\bR_{\ge0}[\overline{E}]+\bR_{\ge0}[\overline{H}_x]+\bR_{\ge0}[\overline{H}_z] & 3a<b\le 4a \\
         \bR_{\ge0}[F]+\bR_{\ge0}[\overline{E}]+\bR_{\ge0}[\overline{H}_z] & b\ge4a
        \end{cases}
    \end{equation*}
    
\end{itemize}
\end{proof}

\begin{pro}
If the weight of $\lambda$ in the \hyperref[unigonal case 1']{Case $1'$}, then $S$-function of the valuation $\ord_F$  is  given by
\begin{equation*}
S_{(X,cC)}(\ord_F)=\frac{106b+83a}{48} (1-2c)
\end{equation*}
\end{pro}

\begin{proof}

By (\ref{scaling}), we need to compute the volume of $L_t=-\mu^*K_X-tF$. To do so, finding the positive part of $L_t$ defined in Proposition \ref{zariski} is sufficient.
Direct computation shows
\begin{equation}\label{intersection number2}
    L_t.F=\frac{t}{ab},\ \ L_t.\overline{E}=1-\frac{t}{b},\ \ L_t.\overline{H}_y=\frac{1}{2}-\frac{t}{a}.
\end{equation} Thus, according to Lemma \ref{Curve Cone},  $L_t$ is nef when $0\le t\le\min\{\frac{a}{2},b\}$. On the other hand, observe $\sqrt{8ab}>\min\{\frac{a}{2},b\}$ and thus \[L_t^2=8-\frac{t^2}{ab}\ge0 \] for $0\le t\le \min\{\frac{a}{2},b\}$.  So we get
 \[ P(L_t)=L_t, \ \ 0\le t\le \min\{\frac{a}{2},b\}  \]
 and need to compute the positive part $P(L_t)$ for $t\ge\min\{\frac{a}{2},b\}$.

 We divide the remaining computations into two cases.
 \begin{itemize}
     \item If $\min\{\frac{a}{2},b\}=b$, from (\ref{intersection number2}) $ L_t$ fails to be nef for $t>b$ since $L_t.\overline{E}<0$ for $t>b$. So we may assume the negative part for $t>b$ is $s \overline{E}$ for some $s>0$ and thus the positive part $P_2(L_t)=L_t-s\overline{E}$.  By Proposition \ref{zariski}, we know that \[P(L_t).\overline{E}=(L_t-s\overline{E}).\overline{E}=0,\] which implies \[
     s=\frac{L_t\overline{E}}{(\overline{E})^2}=\frac{t-b}{a+b}.
     \] Then by direct computation,  the intersection numbers of $P(L_t)$ with $F,\overline{E}$ and $\overline{L}$ are given by
\[ P(L_t).F=\frac{t+a}{a(a+b)},\ \ \ P(L_t).\overline{E}=0,\ \ \ P(L_t).\overline{H}_y=\frac{1}{2}-\frac{t}{a} \]
     which implies that $P(L_t)$ is nef $b\le t\le\frac{a}{2}$. So we have
     \[P(L_t)^2=L_t^2-sL_t\cdot \overline{E} =8-\frac{t^2}{ab}-\frac{t-b}{a+b}(1-\frac{t}{b})> 0,\  \  \ \hbox{for}\  b\le t \le \frac{a}{2}. \]
And we need to continue finding the positive part of $P(L_t)$ for $t>\frac{a}{2}$. The same trick to find the positive part  $L_t$ for $b\le t \le \frac{a}{2}$ will be applied again. We assume the negative part for $P(L_t)$  is 
$s'\cdot \overline{H}_y$ for some $s'>0$ due to the fact  $P(L_t). \overline{H}_y<0$ for $t>\frac{a}{2}$. Then from \[P(L).\overline{H}_y=(P(L_t)-s'\overline{H}_y).\overline{H}_y=0,\] we get $s'=\frac{4t-2a}{3a+4b}$  and 
\[P(L_t)^2=(L_t-s\cdot \overline{E}-s' \cdot \overline{H}_y)^2=\frac{(5a+6b)^2}{(a+b)(3a+4b)} .\]
 Moreover, $P(L_t).F,\ P(L_t).\overline{E},\ P(L_t).\overline{H}_y$ are non-negative for $\frac{a}{2}\le t\le5a+6b$.  Thus,  we conclude that  the pesudo-effective threshold of $L_t$ is \[ \tau= 5b+6a .\] In a summary,  the Zariski decomposition of $L_t$ in this case is given by
      \begin{equation*}
    P(L_t)=
    \begin{cases}
      L_t , & \ \ 0\le t\leq b, \\
    L_t-\frac{t-b}{a+b}\overline{E}  , & \ \ b\le t\le\frac{a}{2},\\
   L_t-\frac{t-b}{a+b}\overline{E}-\frac{4t-2a}{3a+4b}\overline{H}_y, &\ \  \frac{a}{2}\le t\le5a+6b
    \end{cases}
\end{equation*}

\item if $\min\{\frac{a}{2},b\}=\frac{a}{2}$, the computation is parallel to the case $\min\{\frac{a}{2},b\}=b$ and we leave it to the interested reader. In the case, we have Zariski decomposition
     \begin{equation*}
    P(L_t)=
    \begin{cases}
      L_t , & \  0\le t\leq \frac{a}{2}, \\
    L_t-\frac{4t-2a}{3a+4b}\overline{H}_y  , &\  \frac{a}{2}\le t\le b,\\
   L_t-\frac{4t-2a}{3a+4b}\overline{H}_y-\frac{t-b}{a+b}\overline{E}, &\  b\le t\le5a+6b
    \end{cases}
\end{equation*}
 A direct computation shows that in both cases we have
 \[ \int_{0}^{\infty}\vol(L_t) dt=\int_{0}^{\infty}P_t^2 dt=\frac{106a+83b}{6}.  \] So we conclude that \[S_{(\bbF_1,cC)}(F)=\frac{1-2c}{8}\cdot \frac{106a+83b}{6}=\frac{106a+83b}{48}(1-2c).\]
\end{itemize}
\end{proof}

For \hyperref[unigonal case 2']{Case $2'$} and \hyperref[unigonal case 3']{Case $3'$}, the computation is similar.  We omit it and state the computational result as follows
\begin{pro}
The $S$-function of the valuation $\ord_F$ in \hyperref[unigonal case 2']{Case $2'$} is  given by
\begin{equation*}
S_{(X,cC)}(\ord_F)=
  \begin{dcases}
    (18\sqrt{a(a+b)}-\frac{b+26a}{3})\frac{1-2c}{8}, & b<3a;\\
    \frac{25b+83a}{48} (1-2c), & b\ge3a.
    \end{dcases}
\end{equation*}

\end{pro}

\begin{pro}
The $S$-function of the valuation $\ord_F$ in \hyperref[unigonal case 3']{Case $3'$} is given by 
\begin{equation*}
    S_{(X,cC)}(\ord_F) = 
    \begin{dcases}
    \frac{72a+27b+4\sqrt{a(4a-b)}}{48}(1-2c), & b<3a;\\
    \frac{82a+25b}{48}(1-2c), & 3a\leq b <4a;\\
    \frac{2\sqrt{b(b-3a)}+110b +375a}{216}(1-2c) , & b>4a.
    \end{dcases}
\end{equation*}
\end{pro}

\begin{pro}\label{-4curve}
    If $C \subset X=Bl_{[1,0,0]}\bP(1,1,4)$ is a curve in $|-2K_X|$ and $B=\pi(C)$ passing the $\frac{1}{4}(1,1)$ singularity $[0,0,1]$, then $(X,cC)$ is K-unstable for any $c\in (0,\frac{1}{2})$.
\end{pro}
\begin{proof}
    Let $\mu: \widetilde{X} \rightarrow X$ be the minimal resolution of the $\frac{1}{4}(1,1)$ singularity with exceptional divisor $F$. Then 
    \begin{equation*}
        K_{\widetilde{X}}=\mu^\ast K_X-\frac{1}{2} F.
    \end{equation*}
 Denote $H_y$ the proper transform of the line $\{y=0\}$ on $\bP(1,1,4)$  and thus 
 \[ \mu^\ast \pi^\ast \{y=0 \} =H_y+E+\frac{1}{4} F \]
Since the curve $B$ passes though $[0,0,1]$, the $B$ is of the following form 
\[  \{z^2f_4(x,y)+ zf_8(x,y)+f_{12}(x,y)=0\} .\]
 Therefore $\ord_F(C) \ge 1$, in particular, 
 \begin{equation*}
     A_{(X,cC)}(F)=\frac{1}{2}-c\cdot \ord_F(C) \le \frac{1}{2}-c.
 \end{equation*}
On the other hand,  $\NS(\widetilde{X})$ is generated by negative curves  $H_y, E$ and $F$ with 
 \begin{equation*}
    \left(
        \begin{array}{c|c c c }
    & F & E & H_y \\
    \hline
     F  & -4 & 0 & 1 \\
    
   E &  0  & -1 & 1 \\
    
  H_y &  1 & 1 & -1\\
    
    \end{array}
    \right)
    \end{equation*}
    Then we deduce that 
    \begin{equation*}
        \int \vol(-\mu^\ast K_X-tF)=\int_0^{\sqrt{2}} (8-4t^2)dt=\frac{16\sqrt{2}}{3}.
    \end{equation*}
Thus by the formula (\ref{scaling})
\[ S_{(X,cC)}(F)=\frac{2\sqrt{2}}{3}(1-2c)> A_{(X,cC)}(F).\]
This proves the pair $(X,cC)$ is destabilised by the valuation $\ord_F$ according to Fujita-Li's criterion. 
\end{proof}

\subsubsection{S-function on index $3$ del Pezzo pairs and their stability}
Recall from the Section \ref{index3surface} , we know $X$ has a unique quotient singularity of type $\frac{1}{9}(1,2)$. We take a weighted blowup $\phi: Z \rightarrow X$ of weight $(1,2)$ with exceptional divisor $F$ and then $Z$ has a $A_1$ singularity in $F$, then we continue to take the blowup at $A_1$ singularity and then the minimal resolution $\pi: M \rightarrow X$ is a composition of $\phi$  and blowup $q:M\rightarrow Z$ at the $A_1$ singularity .

\begin{pro} \label{index3}
Let $(X,C)$ be a index $3$ del Pezzo pair of degree $8$, then it is K-unstable for all $c<\frac{1}{2}$.
\end{pro}
\begin{proof}
We are going to show  $\beta(F)=A_{(X.cC)}(F)-S_{(X.cC)}(F) <0$ for all index 3 del Pezzo pairs. By the proof in Proposition \ref{Cartier index} and formula (\ref{wbu}), 
\begin{equation*}
   A_{(X.cC)}(F)=A_X(F)-c\cdot \ord_F(C)=\frac{1}{3}- \frac{2}{3}\cdot c.
\end{equation*}
Observe that $\vol(-\phi^\ast K_X-t\cdot F)=\vol(q^\ast(-\phi^\ast K_X-t\cdot F))$ and 
\[L_t:=q^\ast(-\phi^\ast K_X-t\cdot F)=-\pi^\ast K_X-\frac{2}{3}F_1-\frac{1}{3}F_2-t\cdot q^*F. \]
Therefore, we will calculate the volume function $\vol(L_t)$ for $L_t$ on $M$. As
\[ q^*F=F_1+\frac{1}{2}F_2,\]
we get 
\[
L_t=(\frac{4}{3}-t)F_1+(\frac{20}{3}-\frac{t}{2})F_2+6E_1+6E_2
\]
By the intersection matrix on $M$ mentioned in  (\ref{M}),  $L_t$ is nef for $t\ge0$. Thus we have $P(L_t)=L_t$ for $t\ge0$ and $\vol(L_t)=(L_t)^2=8-\frac{9}{2}t^2$. Hence  the psedudo effective threold of $L_t$ is $\frac{4}{3}$. So we get 
\begin{equation*}
    \begin{split}
      S_{(X,cC)}(F)= \frac{(1-2c)}{8} \int_0^{\frac{4}{3}} (8-\frac{9}{2}t^2)dt 
      = \frac{8}{9}(1-2c)>A_{(X,cC)}(F)=\frac{1}{3}- \frac{2}{3}\cdot c.
    \end{split}
\end{equation*}
This  finishes the proof.

\end{proof}

\subsection{The K-moduli walls}
We apply the computation of $S$ function to determine all the walls for the K-moduli $\overline{P}_c^K$.

\begin{pro}\label{ckw1}
Let $C \in |-2K_{\bbF_1}|$  be the curve on $\bbF_1$ such that $\pi(C)$ is one of the curves on $\bP^2$ given in the Table \ref{Kwall1}, then the stability threshold of $(\bbF_1,C)$ is a point in $[0,\frac{1}{2}]$ given by 
 \[\{\ \frac{1}{14},\frac{5}{58},\frac{1}{10},\frac{7}{62},\frac{1}{8},\frac{5}{34},\frac{1}{6},\frac{7}{38},\frac{1}{5},\frac{5}{22} ,\frac{2}{7}\ \}. \] In particular, these numbers are K-moduli walls.
\end{pro}
\begin{proof}
We check the curves given in the Table \ref{Kwall1} case by case. We first give a proof for the first critical value $c=\frac{1}{14}$. 
Note that for the curve $C\subset Bl_p \bP^2$, whose image in $\bP^2$  is the plane curve $\{z^4xy=0\}$, we have 
\[C=H_x+H_y+4H_z.\]
By Fujita-Li's criterion \ref{FL} and Proposition \ref{sfunction1},  
\begin{equation*}
    A_{(\bbF_1,cC)}(H_z)=1-4c \ge S_{(\bbF_1,cC)}(H_z)=\frac{5}{6}(1-2c).
\end{equation*}
Thus, we get $ c \le \frac{1}{14}$.
Similarly, for the divisorial valuation given by $H_x$, we have 
\begin{equation*}
  A_{(\bbF_1,cC)}(H_x)=1-c \ge S_{(\bbF_1,cC)}(H_z)=\frac{13}{12}(1-2c)
\end{equation*}
and thus we get $ c \ge \frac{1}{14}$. This shows the stabilty threshold of  $(\bbF_1,C)$ is either empty or $\{ \frac{1}{14} \}$. Thus, it remains to show the pair $(\bbF_1,\frac{1}{14}C)$ is K-semistable.  It is sufficient to show $\beta_{(\bbF_1,\frac{1}{14}C)}(F)\ge 0$ for any plt type blowup divisor $F$ over $\bbF_1$ by Fujita's criterion.  It is well known that \ any plt type blowup is weighted blowup.  Then by  Proposition \ref{S1} and Proposition \ref{S2}, one can easily cheack that
\begin{equation*}
    \begin{split}
        \beta_{(\bbF_1,\frac{1}{14}C)}(F)=a+b - S_{(\bbF_1,\frac{1}{14}C)}(\ord_F) \ge 0
    \end{split}
\end{equation*}

Observe that for other $w$ the pair  $(\bbF_1,C)$  given by  Table \ref{Kwall1} is of complexity one. We will  apply equivariant K-stability criterion in theorem \ref{thmcomplexity1} to show the stability threshold for  $(\bbF_1,C)$  is exactly $\{ w\}$. Let us show for $w=\frac{5}{58}$  and leave the remaining cases to interested reads.  Denote $\lambda=(1,t^2,t^3)$ the 1-PS  acting on $\bP^2$, which can be lifted to $\bbF_1$, then the
Futaki character is just the $\beta$-invariant of divisor valuation $\ord_F$.  Here $F$ is the exceptional divisor of weighted blowup of $\bbF_1$ at the  point $([1,0,0],[1,0])$ under the coordinate in Case \ref{Case2}. In this case, the  weight for the blowup is  $a=2,\ b=1$, then by Proposition \ref{S1}, 
\begin{equation*}
    \begin{split}
       \beta_{(\bbF_1,\frac{5}{58}C)}(\ord_F)& = A_{(\bbF_1,\frac{5}{58}C)}(\ord_F)- S_{(\bbF_1,\frac{5}{58}C)}(ord_F) \\
      &  = (a+b)-2\cdot \frac{5}{58}-\frac{14a+13b}{12}(1-2\cdot \frac{5}{58})=0
    \end{split}
\end{equation*}
 The $\lambda$-vertical divisors $F$ are given by  $H_x,H_y,H_z$ and $E$,  $L=\{ y-z=0\}$, then by Proposition \ref{sfunction1}, it is easy to see $\beta_{(\bbF_1,\frac{5}{58}C)}(\ord_F)>0$. This finishes the proof of K-polystablity of $(\bbF_1,\frac{1}{14}C)$ by theorem \ref{thmcomplexity1}.

\end{proof}
\begin{pro}\label{ckw2}
Let $C \in |-2K_{Bl_p\bP(1,1,4)}|$  be the curve on $Bl_p\bP(1,1,4)$ such that $\pi(C)$ is one of the curve on $\bP(1,1,4)$ given in the Table \ref{Kwall2}, then the stability threshold of $(Bl_p\bP(1,1,4),C)$ is a point in $(0,\frac{1}{2})$ given by 
 \[\{\ \frac{29}{106},\frac{31}{110},\frac{2}{7}, \frac{35}{118} \ \}. \] In particular, these numbers are K-moduli walls.
\end{pro}
\begin{proof}
The computation to check the curve on $\bP(1,1,4)$ given in the Table \ref{Kwall2} is similar to that in Proposition  \ref{ckw1}. We omit the computation. 
\end{proof}

\subsection{Proof of part (1) of theorem \ref{mainthm1}} By the local VGIT structure of K-moduli space $\overline{P}_c^K$ (see theorem \ref{kwc}), we know each log del Pezzo pair $(X,C)$ parametrized by the center of each wall $w$ is a $\bG_m$-equivariant degeneration of  the log del Pezzo pairs on exceptional loci. Thus,  $(X,C)$ admits a $\bG_m$-action. Moreover, by Proposition \ref{classification of degeneration} and Proposition \ref{index3},  such $X$ is either $\bbF_1$ or $Bl_p\bP(1,1,4)$. Therefore to prove part (1) of theorem \ref{mainthm1}, it remains to show the curves listed in the Table \ref{Kwall1} and Table \ref{Kwall2} are all curves such that the pair $(X,C)$ admits  $\bG_m$-action and its K-stability threshold is a point. This is done by Proposition \ref{ckw1}, Proposition \ref{ckw2} and the following algorithm \footnote{We write an easy Python code to help us find the wall. The code for walls on $\bbF_1$ can be found  \href{https://changfeng1992.github.io/SiFei/wall test - (z,u) - F1.ipynb}{here} and for $Bl_p\bP(1,1,4)$ is \href{https://changfeng1992.github.io/SiFei/wall test - (z,u) - Bl P(1,1,14).ipynb}{here}.}  to find the potential critical curves.  Let $\lambda$ be the $\bG_m$-action with weight $(\lambda_1,\lambda_2,\lambda_3)$ as in Section \ref{4.1.1} and Section \ref{4.1.2}. If the pair $(X, C)$ is in the center of wall $w$,  then the $\beta$-invariant of the valuation associated to  $\lambda$ is zero. In particular,
\begin{equation}\label{algorithm}
    A_{(X,w\cdot C)}=S_{(X,w\cdot C)}.
\end{equation}
Now we use the equation (\ref{algorithm}) to give algorithm to find the potential walls for $X=\bbF_1$. The same algorithm works for $X=Bl_p \bP(1,1,4)$ and we omit the details.  As in the same coordinate chart of the computation of $S$-function in Proposition \ref{S1} and Proposition \ref{S2},    each curve $C \in |-2K_{\bbF_1}|$  is of the form
     \[ C= \pi^{-1}( D) -2E\] where 
    $D$ is a curve on $\bP^2$ defined by
    \begin{equation}\label{planesextic}
        x^4f_2(y,z)+x^3f_3(y,z)+x^2f_4(y,z)+xf_5(y,z) +f_6(y,z)=0.
    \end{equation}
    If $D$ is defined by one monomial in (\ref{planesextic}), say $\{y^iz^jx^{6-i-j}=0\}$, then  it is easy to check that  $\{x^4yz=0\}$ is the only possibility. Indeed, by Proposition \ref{sfunction1}, we  have
    \[
    \beta_{(\bbF_1,wC)}(H_x)=A_{(X,wC)}-S_{(X,wC)}=1-w\cdot i-\frac{13}{12}(1-2w)\ge0
    \]
    Thus $0\le i\le 1$ since $w<\frac{1}{2}$ and the same reason imply for  $0\le j \le 1$. So the unique solution is $i=j=1$.
   
    Otherwise, $D$ is defined by more than two  monomials.  By the  $\bG_m$- invariance on the curve $C$,   any two monomials in the defining equation of $D$, say $y^iz^jx^{6-i-j}$ and  $y^{i'}z^{j'}x^{6-i'-j'}$ have the same weights under $\lambda$, that is, 
    \begin{equation}\label{monoialequ}
        i \cdot \lambda_2+j\cdot \lambda_3+(6-i-j)\cdot \lambda_1= i' \cdot \lambda_2+j'\cdot \lambda_3+(6-i'-j')\cdot \lambda_1
    \end{equation}
Then we divide the remaining discussion  according to the different chart in the computation of $S$-function of divisorial valuation.
\subsubsection*{ In the Case (\ref{Case1}).}
Recall that $\ord_F$ is given by the exceptional divisor of the weighted blowup. For simplicity wo only consider the weighted blowup along the point $[0,1,0]$ with weight
    \[\wt(x)=\lambda_1-\lambda_2=a,\ \ \ \wt(z)=\lambda_3-\lambda_2=b\] 
By the equation (\ref{monoialequ}), we may set $a=(j'-j), b=i'+j'-(i+j)$.
We have
\[ A_{(\bbF_1, wC)}=a+b-m\cdot w \]
       where $m=(6-i-j)a+jb$.
       By Proposition \ref{S1} and 
        equation (\ref{algorithm}), we get 
    \begin{large}
       \begin{equation*}
            w=\begin{cases}
          \frac{b^2}{12am-24a^2-24ab-2b^2}, &0<b<a\\
          \frac{2b-a}{12m-26a-20b}, & b\ge a
        \end{cases}
        \end{equation*} 
       \end{large} 
       
      Let $y^iz^jx^{6-i-j}$ and  $y^{i'}z^{j'}x^{6-i'-j'}$ go through all monimals appeared in equation (\ref{planesextic}), then we obtain all potential walls $w$ and corresponding curves $C$ in this case.
      
      \subsubsection*{ In the Case (\ref{Case2}).}
      Recall that in this case, then $F$ is the exceptional divisor of  weighted blowup $Y\rightarrow\bbF_1$ where the center is $([1,0,0],[1,0])$ and  the local parameter is $(z,u)$ with weight $(\lambda_3-\lambda_1,\lambda_2-\lambda_3)$. (We can assume $\lambda_2>\lambda_3$. If $\lambda_3<\lambda_2$, the computation is similar). We first write the equation of $C$ in terms of $(z,u)$ as follows:
      
      \begin{equation}\label{strict transform}
       \widetilde{f}_2(u)+z\widetilde{f}_3(u)+z^2\widetilde{f}_4(u)+z^3\widetilde{f}_5(u)+z^4\widetilde{f}_6(u)=0   
      \end{equation}
where $\widetilde{f}_k(u)=z^{-k}f_k(zu,z)$.  
      Then $y^iz^jx^{6-i-j}$ and  $y^{i'}z^{j'}x^{6-i'-j'}$ should be replaced by $z^{i+j-2}u^i$ and $z^{i'+j'-2}u^{i'}$. They have same weight. Let $a=\lambda_3-\lambda_1$, $b=\lambda_2-\lambda_3$. We have
      \begin{equation*}\label{a&b}
          (i+j-2)a+jb=(i'+j'-2)a+j'b.
      \end{equation*}
      So we may set  $a=j'-j$ and $b=i+j-(i'+j')$. As 
      \[
      A_{(\bbF_1, wC)}=a+b-m\cdot w,\ m=(i+j-2)a+jb
      \]
    by Proposition \ref{S2}, we get
     \begin{equation*}\label{wall}
         w=\frac{2a+b}{28a+26b-12m}.
     \end{equation*}
 As before  let  $z^{i+j-2}u^i$ and $z^{i'+j'-2}u^{i'}$ go through all monomials appeared in equation (\ref{strict transform}), we obtain all potential walls $w$ and the corresponding curves $C$ in this case.

\begin{rem}
  For the above algorithm in the case $X=Bl_{[1,0,0]}\bP(1,1,4)$, we  only need to consider the local equations for curve $C$ containing at least two monomials.  Since by Proposition \ref{-4curve}, the monomial $z^3$ for $B$ must appear and it is also easy to check  $(X,cC)$ is unstable for any $c\in (0,\frac{1}{2})$ if  $B=\{z^3=0\}$ by the divisorial valuation $\ord_{H_z}$.
\end{rem}
\begin{rem}
   For the last wall $w=\frac{35}{118}$,  the general equation for $B$ obtained by the algorithm is 
   \[ z^3+a_1z^2yx^3+a_2zy^2x^6+a_3y^3x^6=0.\]
   But under the transform $z\mapsto z+a yx^3$ for suitable $a\in \bC$,  the equation will be equivalent the one in the Table \ref{Kwall2}. 
\end{rem}
\vspace{0.2cm}
\section{Explicit wall-crossings for K-moduli of degree 8 log Fano pairs} \label{explicitwc}
In this section, we will describe the wall-crossings for  K-moduli space $\overline{P}^K_c$ explicitly. 
\subsection{Wall-crossings on  surface $\bbF_1$}

\subsubsection{ First wall $w=\frac{1}{14}$}\label{section1stwall}
In the  case   del Pezzo pairs of deg $8$  where the surface is  $X\cong \bP^1 \times \bP^1$,   the K-moduli space are  studied in \cite{ADL2021} and is shown always nonempty for  $c>0$.  Our case is a little  different and we first show 
\begin{pro}\label{1}
For $0<c<c_0=\frac{1}{14}$, then $\overline{P}^K_c$ is empty.
\end{pro}
\begin{proof}
To prove the first Proposition, it is enough to show that for any $c<\frac{1}{14}$ there exists a divisor $F$ over $(\bbF_1, cC)$ such that $\delta_{(\bbF_1, cC)}(F)<1$. We use the notation as (\ref{Case2}). Let $\mu\colon Y\rightarrow\bbF_1$ be a weighted blowup of $\bbF_1$ which centered at $([1,0,0],[1,0])$ and the weight on $z$ and $u$ are $a=1$ and $b=2$ respectively. From the formula (\ref{wbu})
\[
    A_{(\bbF_1,cC)}(F)=3-c(i_0+2 j_0),
\]
By Proposition \ref{S2}, we know that $S_{(\bbF_1,cC)}(F)=\frac{10}{3}(1-2c)$. 
The Fujita-Li's criteria implies if $((\bbF_1,cC))$ is K-semistable, then
\[
3-c(i_0+2j_0)\ge\frac{10}{3}(1-2c).
\]
It turns out that $$c \ge \frac{1}{20-3i_0-6j_0}$$ Note that the local equation for $C$ at $([1,0,0],[1,0])$ in the coordinate $(z,u)$  is \[f_2(1,u)+zf_3(1,u)+\cdots+z^4f_6(1,u) =0 .\] Thus  $i_0 \ge 0 $ and $j_0 \ge 1$ which implies $c\ge\frac{1}{14}$. This finishes the proof. 
\end{proof}
Before description of the  second wall,  let us identify the K-moduli space as a GIT quotient space $\bP V \q \bC^\ast$. First we give the construction for  $\bP V \q \bC^\ast$.  Recall $\pi: \bbF_1 \rightarrow \bP^2$ is the blowup at $p=[1,0,0]$ and $\Delta \subset |\cO_{\bP^2}(6)|$ be the discriminant locus of plane sextic curves, then there is a universal nodal curve with two projections  $p_1: \cv \rightarrow \bP^2 $ and $p_2: \cv \rightarrow \Delta$.  So $p_2(\cv_p)$ is the locus of plane sextic curves singular at $p$ where $\cv_p:=p_1^{-1}([1,0,0])$. It is not hard to see $$p_2(\cv_p) \cong |\cO_{\bP^2}(6) \otimes \mathfrak{m}_p^2| \cong \bP^{24}.$$
Note that $|\cO_{\bP^2}(6) \otimes \mathfrak{m}_p^2|$ is  also the  parameter space of the polynomials of the form
\begin{equation}\label{nodalequ}
    x^4f_2(y,z)+x^3f_3(y,z)+\cdots +f_6(y,z).
\end{equation}
 The  non-reductive group $\aut(\bbF_1)=\GL(2,\bC)  \rtimes \bC^2$ acts on $|\cO_{\bP^2}(6) \otimes \mathfrak{m}_p^2|$ by
\begin{equation*}
    \begin{split}
        y \mapsto & a_{11}y+a_{12}z \\
        z \mapsto & a_{21}y+a_{22}z\\
        x\mapsto &  x+b_1y+b_2z, (b_1,b_2)\in \bC^2.
    \end{split}
\end{equation*}
Let $U\subset |\cO_{\bP^2}(6) \otimes \mathfrak{m}_p^2|$ be the locus where $f_2(y,z)$ is a smooth conic, i.e., the sextic curve with a node at $p$, then modulo  the action of $\aut(\bbF_1)$, each $f\in U$ has the  normal form $f=x^4yz+h$ where 
\begin{equation}\label{nf}
\begin{split}
     &h= x^3\widetilde{f}_3(y,z)+x^2f_4(y,z)+xf_5(y,z)+f_6(y,z),\\  &\widetilde{f}_3(y,z)=a_1yz^2+a_2y^2z
\end{split}
\end{equation}
Denote the $V$ by the $ \bC$-vector space spanned by the monomials in the normal form.
Thus, $U$ is the orbit $\aut(\bbF_1)\cdot  \bP W$.
Moreover, the stabilizer of $\bP W$ is a 2-dimensional torus $T = \left\{\mdiag(a,b,c)\,\middle| \,  abc =1 \right\}\cong \bC^* \times \bC^\ast$.
Therefore, we have a reductive GIT space \begin{equation} \label{gitd8}
    \bP V \q T.
\end{equation}

\begin{thm}\label{2nd}
There is an isomorphism \[
    \overline{P}^K_c \cong \bP V \q \bC^\ast 
\]  
 for $\frac{1}{14}<c<\frac{5}{58}$.  
\end{thm}
\begin{proof}
To show the K-moduli is  isomorphic to a GIT quotient space, we follow the arguments in \cite[proof of therorem 1.1]{LiuX}. First, we claim that for any K-semistable pair $(X,cC)$ in $\overline{P}^K_c$ for  $\frac{1}{14}<c<\frac{5}{58}$, we have \[ X\cong Bl_{[1,0,0]} \bP^2,\  \ \  \pi(C)= \{f=0\} .\] where $f\in V-\{ 0 \}$.  Indeed, by Corollary \ref{classification of degeneration} and Proposition \ref{index3}, $X$ is either $\bbF_1$, $Bl_p(\bP(1,1,4))$. If $X=Bl_p\bP(1,1,4)$, by taking  divisorial valuation $\ord_{H_y}$, we get 
\[ c \ge \frac{29}{106-24 \cdot \ord_{H_y}(C)  } > \frac{5}{58},\] which  is impossible.  This shows $X\cong Bl_{[1,0,0]} \bP^2$. If 
\[\pi(C):=\{ x^4f_2(y,z)+x^3f_3(y,z)+\cdots+f_6(y,z)=0\} \subset \bP^2\] is not nodal at $[1,0,0]$, that is the quadratic $f_2(x,y)$ has rank $\le 1$. Then we take a weighted blowup of $\bbF_1$ at $[1,0,0] \times [1,0]$ with weight $(2,1)$ under local coordinate $(y,u)$ and $F$ is the exceptional divisor.  By Proposition \ref{S2}, 
\[ \beta_{(\bbF_1,cC)}(F) \le 3- 2c - \frac{41}{12}(1-2c)=\frac{29}{6}c-\frac{5}{12}<0\]
for $c<\frac{5}{58}$. This finishes the proof of the claim. 

Next, there is a universal family of surface pairs $(\cX,\cC)$ over $\bP V$ where $\cX=\bbF_1 \times \bP V$ and  $\cC \subset \bbF_1 \times \bP V$ is the restriction of  hypersurfaces 
\[  \{x^4yz+x^3f_3(y,z)+...+f_6(y.z)=0\}  \subset  \bbF_1 \times \bP (V_2 \oplus \cdots \oplus V_6) \]
on $ \bbF_1 \times \bP V$ under the closed embedding $ \bbF_1 \times \bP V \hookrightarrow \bbF_1 \times  \bP (V_2 \oplus \cdots \oplus V_6)$. By the claim, the parameter space  $\cP_c$ of K-semistable del Pezzo pairs can be realised as a subset of $\bP V$. By direct computation,  the CM line bundle $\lambda_{CM,c}$ is propositional to  the restriction of $\cO_{\bP V}(1)$. Thus, we get injective morphism 
\begin{equation}\label{parameter}
  \cP_c \ \hookrightarrow \  (\bP V)^{ss} 
\end{equation} 
By properness of K-moduli stack (see \cite[Theorem 3.1]{ADL19} for smoothable case and \cite{LXZ} for general case),  the morphism (\ref{parameter}) descends to an isomorphism of  $  \overline{P}^K_c \cong \bP V \q \bC^\ast
$. 
\end{proof}
\begin{rem}
The reductive GIT quotient space $\bP V \q \bC^\ast $  is  isomorphic to  the non-reductive GIT quotient space $\lvert-2K_{\bbF_1}\rvert \q \aut(\bbF_1)$  in the sense of Doran-Kirwan \cite[definition]{DK}.
\end{rem}

    

\subsubsection{The second wall $w=\frac{5}{58}$ and the 1st divisorial contraction}

Denote $C^{-}_1$ and $C^0_1$ the curve on $\bbF_1$  such that   $\pi(C^{-}_1)=\{x^4z^2+x^4zy+x^3y^3=0\}$ and $\pi(C^{0}_1)=\{x^4z^2+x^3y^3=0\}$ respectively.

\begin{pro}\label{1stdivc}
At the wall $w=\frac{5}{58}$, there are natural  birational morphisms 
\[ \overline{P}^K_{\frac{5}{58}-\epsilon}\  \xlongrightarrow{p^-}\  \overline{P}^K_{\frac{5}{58}}  \  \xlongleftarrow{p^+}  \  \overline{P}^K_{\frac{5}{58}+\epsilon}\]
where 
\begin{enumerate}
    \item $p^-$ is an isomorphism. More precisely, $p^-$ is identity outside the point representing K-polystable pair $(\bbF_1, \frac{5}{58}C^{0}_1)$ and maps K-polystable  pair $(\bbF_1, (\frac{5}{58}-\epsilon)C^{-}_1)$ to   $(\bbF_1, \frac{5}{58}C^{0}_1)$ .
    \item $p^+$ is a Kirwan type blowup at the point $[(\bbF_1, \frac{5}{58}C^{0}_1)]$.  The  exceptional divisor $E_c^+$ parametrizes $S$-equivalence classes of K-semistable  pairs $(\bbF_1,C)$ where $B=\pi(C)$ is described in the second row  in  the table  \ref{tabAE}.  Moreover, $E_c^+$ is birational to the hyperelliptic divisor  $\NL(A_2) \cong H_h$ in $\cF$.
\end{enumerate}
\end{pro}
\begin{proof}
We first show that $p^-$ is an isomorphism. It is clear that $p^-$ is a birational morphism between normal proper varieties since there exists a common open subset shared by $\overline{P}^K_{\frac{5}{58}}$ and $\overline{P}^K_{\frac{5}{58}-\epsilon}\cong\overline{P}^{GIT}$. Indeed, the Picard number $\rho(\overline{P}^{GIT})$ is one. 
Note that we have the following relations:
\[
\Pic(\bP V \q \bC^\ast)\hookrightarrow\Pic_{\bC^*}(\bP V^{ss})\hookrightarrow\Pic(\bP V^{ss})
\]
Since we have a surjective morphism from $\Pic(\bP V)\cong\bZ$ to $\Pic(\bP V^{ss})$, we know that $\rho(\overline{P}^{GIT})=1$. Thus $p^-$ is isomorphism between good moduli spaces.   Now we give the explicit description of K-semistable replacement for $p^-$ and $p^+$.
Note that as in the proof of theorem \ref{2nd},  the pair $(\bbF_1, \frac{5}{58}C)$ is K-semistable only if $B=\pi(C)$ is a plane curve with at worst $A_2$-singularity at $p=[1,0,0]$ . Assume $B$ has $A_2$-singularity at $p=[1,0,0]$, then up to a coordinate changes it is a plane curve defined by the following equation
\begin{equation}\label{A2}
    x^4z^2+ x^3f_3(y,z)+ x^2f_4(y,z)+ xf_5(y,z)+f_6(y,z)   =0
\end{equation}
By checking the stability threshold of pairs $(\bbF_1,C)$ via the computation of $S$-function in Proposition \ref{S2}, one can see that $(\bbF_1, C^{0}_1)$ is the only curve with the stability threshold $\{ \frac{5}{58} \}$ and $(\bbF_1, \frac{5}{58}C^{0}_1)$ is K-polystable.
Therefore, outside the point  $[(\bbF_1, \frac{5}{58}C^{0}_1)]\in \overline{P}^K_{\frac{5}{58}-\epsilon}$, the K-moduli stacks for  $\overline{P}^K_{\frac{5}{58}-\epsilon}$ and $\overline{P}^K_{\frac{5}{58}}$  will share the same parameter space and then the morphism of K-moduli stack  will descend to the isomorphism $p^-$ over $\overline{P}^K_{\frac{5}{58}-\epsilon}-[(\bbF_1, \frac{5}{58}C^{0}_1)]$.  Note that under 1-PS $\lambda$ with weight given by $(0,2,3)$,  
\[ \mathop{\lim} \limits_{t \rightarrow 0} \ \lambda(t) C_1^-=C_1^0 \]
This will induces a $\bG_m$-equivariant degeneration $(\fX,\cC) \rightarrow \bA^1$ with generic fiber isomorphic to $(\bbF_1,C^{-}_1)$ and central fiber isomorphic to $(\bbF_1,C^{0}_1)$. 
Thus by the local VGIT interpretation of K-moduli at the point $[(\bbF_1, \frac{5}{58}C^{0}_1)]$ (see \cite[Theorem 3.33]{ADL19}, \cite{DH}), $p^-$ will map the point $[(\bbF_1, (\frac{5}{58}-\epsilon)C^{-}_1)]$ to  the point  $[(\bbF_1, \frac{5}{58}C^{0}_1)]$ and the exceptional locus $E_w^+$ of $p^+$ parametrizes pairs $(\bbF_1,C)$ such that  $\pi(C)$ is a plane curve defined by equation (\ref{A2}) and $ \mathop{\lim} \limits_{t \rightarrow 0} \ \lambda^{-1}(t) C=C_1^0$. 

By the computation of  Néron-Severi group of  K3 surface $(Y,\tau)$ with involution $\tau$ obtained by a generic pair $(\bbF_1,C)$ in $E_w^+$ (see \cite[Section 4.1]{PSW1}),  we know  the isomorphic class $[(Y,\tau)]\in \NL(A_2)$ and so this shows $E_w^+$ is birational to $\NL(A_2)$.  


\end{proof}

\subsubsection{The remaining walls on $\bbF_1$}
\begin{pro}\label{1stflip}
    At the walls $w\in \{  \frac{1}{10},\frac{7}{62},\frac{1}{8},\frac{5}{34},\frac{1}{6},\frac{7}{38},\frac{1}{5},\frac{5}{22} , \frac{2}{7}\} $, there are flips  
\[ \overline{P}^K_{w+\epsilon}\  \xrightarrow{p^+}\  \overline{P}^K_{w}  \  \xleftarrow{p^{-1}}  \  \overline{P}^K_{w-\epsilon}\]
where the center $Z_w$  is either  points or  a rational curve   described in the Table\ref{Kwall1}. Moreover,
\begin{itemize}
    \item $E_w^{-}=(p^{-})^{-1}(Z_w)$ parametrizes $S$-equivalence classes of K-semistable  pairs $(\bbF_1,(w-\epsilon)\cdot C)$  where $B=\pi(C)\subset \bP^2$  is a plane curve of the form  listed in the Table\ref{F1-}.
 \item $E_w^{+}=(p^+)^{-1}(Z_w)$ parametrizes  $S$-equivalence classes of K-semistable  pairs $(\bbF_1,(w+\epsilon)\cdot C)$  where $B$ is listed in the Table\ref{tabAE}, \ref{tabDE} ,\ref{tabE}. $E_w^+$ is birational to a Noether-Lefschetz locus $\NL(L)$ where $L$ is the complement of $(E_7\oplus A_n)$, $(E_7\oplus D_n)$, $(E_7\oplus E_n)$ or the their modifications in $\II$. 
\end{itemize}
  When $c$ changes from $w-\epsilon$ to $w+\epsilon$,  the K-polystable pairs in $E_w^{-}$ is replaced by the K-polystable pairs in $E_w^{+}$.

  \begin{center}
  \renewcommand*{\arraystretch}{1.2}
\begin{table}[ht]
    \centering
      \begin{tabular}{ |c|c|c|}
    \hline
      $c$ &   equation for  curve  $B$ on  $\bP^2$  &   $\dim E^-$ \\ \hline 
      $ \frac{1}{10}-\epsilon$  & $x^4z^2+x^3zy^2+a\cdot x^2y^4+ x^3y (xl(y,z)+by^2)=0$  & {$2$} \\ \hline 
      
      $\frac{7}{62}-\epsilon$ & $x^4z^2+xy^5+x
      ^2y(x^2l_1(z,y)+x y l_2(y,z)+by^3)=0$  & $2$  \\   \hline 

      \multirow{3}{*}{$\frac{1}{8}-\epsilon$}
       &  $x^4z^2+x^2zy^3+a\cdot y^6$ & \multirow{2}{*}{$5$}  \\ 
       & $+xy(x^3 l_1(y,z)+x^2y l_2(y,z)+y^3l_3(x,y))=0$ & \\
       \cline{2-3}
        &  $ x^3f_3(y,z) + x^4f_2(y,z)$  & {$2$} 
        \\ \hline
      
      \multirow{3}{*}{$\frac{5}{34}-\epsilon$}  & $x^3z^2y+x^2y^4+x^4f_2(y,z)+x^3y^2l_1(y.z)=0$ &  $3$ \\  \cline{2-3}
      &  $x^4z^2+xzy^4+x^4yl_1(y,z)+ $   & \multirow{2}{*}{$6$}   \\ 
      & ${x^3y^2 l_2(y,z)}+x^2y^3l_3(y,z)+y^5l_4(x,y)=0$ & \\
      \hline

      \multirow{3}{*}{ $\frac{1}{6}-\epsilon$ } & $x^3z^2y+x^2zy^3+a\cdot xy^5+x^4f_2(y,z)+x^2y^2(l_1(y,z)+by^2)=0$ &  $5$ \\  \cline{2-3} 
      & $x^4z^2+zy^5+$  &  \multirow{2}{*}{$8$}  \\ 
      & $x^4yl_1(y,z)+y^2(x^3l_2(y,z)+x^2l_3(y,z)+xy^2l_4(y,z)+by^4)=0$ & \\
      \hline
      
    \multirow{2}{*}{ $\frac{7}{38}-\epsilon$}  & $x^3z^2y+y^6+x^4f_2(y,z)+x^3y^2l_1(y,z)+x^2y^3l_2(y,z)+xy^5=0$ & $6$  \\  \cline{2-3} 
      &$x^3z^3+x^2y^4+x^4f_2(y,z)+x^3yg_2(y,z)=0$ & $4$\\ \hline
      
     \multirow{2}{*}{$\frac{1}{5}-\epsilon$} &  $x^3z^2y+xzy^4+x^4f_2(y,z)+x^3y^2l_1(y,z)$ &  \multirow{2}{*}{$7$} \\   
     &  $+x^2y^3l_2(y,z)+y^5l_3(x,y)=0$ & \\
     \hline
     
    \multirow{3}{*}{$\frac{5}{22}-\epsilon$}  & $x^3z^3+x^2zy^3+x^4f_2(y,z)+x^3yf_2(y,z)+x^2y^4=0$  &  $5$\\  \cline{2-3} 
    & 
    $x^3z^2y+z y^5+x^3(x f_2(y,z)+y^2 l_1(y,z))$ & \multirow{2}{*}{$9$} \\
     & $+y^3(x^2l_2(y,z)+xyl_3(y,z)+ay^3)=0$ &  \\       
       \hline

      $\frac{2}{7}-\epsilon$ & $x^3z^3+xy^5+x^4f_2(y,z)+x^3yf_2(y,z)+x^2y^3l_1(y,z)=0$ &  $6$  \\   \hline 
       
   \end{tabular}
    \caption{Equations for  curves in  exceptional locus $E^-$}
    \label{F1-}
\end{table}
\end{center}

\end{pro}
\begin{proof}

The proof is similar to Proposition \ref{1stdivc}. Let \[ w\in \{  \frac{1}{10},\frac{7}{62},\frac{1}{8},\frac{5}{34},\frac{1}{6},\frac{7}{38},\frac{1}{5},\frac{5}{22} , \frac{2}{7}\} .\] be the walls. By part (1) of theorem \ref{mainthm1}, the center $Z_w$ of  $\overline{P}^K_w$ is known in Table \ref{Kwall1}, which is either a point or a rational curve.  Let $C_0$ be a curve listed by the Table \ref{Kwall1} at wall $w$   Then  using the local VGIT presentation of   K-moduli space $\overline{P}^K_w$ at the point $[(\bbF_1,wC_0)]$ in center $Z_w$,  we get the description of local behavior at $[(\bbF_1,wC_0)]$ as follows:  $E_w^-$ parametrizes the curves $C^-$ on $\bbF_1$ such that the limits  $\mathop{\lim} \limits_{t \rightarrow 0} \ \lambda(t) C^-=C^0$ under 1-PS $\lambda$ of the corresponding weight in the Table\ref{Kwall1}. In particular, we get the equation of curve $C^-$ listed in Table \ref{F1-}.  One can use such 1-PS  $\lambda$ to produce special degeneration from the pair $(\bbF_1,C^-)$ to the pair $(\bbF_1,wC_0)$. 
$E_w^+$ is described similarly  by using 1-PS $\lambda^{-1}$, that is,   $E_w^-$ parametrizes the curves $C^+$ on $\bbF_1$ such that the limits  $\mathop{\lim} \limits_{t \rightarrow 0} \ \lambda^-(t) C^+=C^0$. Note that the stabilizer of the pair  $(\bbF_1,C)$ is a one dimensional torus.  Then by the following  dimension formula \cite[Theorem 0.2.5]{DH} in VGIT
\begin{equation}\label{dimvgit}
    \dim E_w^-+\dim E_w^+=17+\dim Z_w,
\end{equation}
  the dimension for $E^-$ is computed. In this way we finish the proof.
\end{proof}

\subsection{Wall-crossing on surface $Bl_p\bP(1,1,4)$}

\subsubsection{Wall  $w=\frac{29}{106}$ and the 2nd  divisorial contraction}

We will construct an explicit 
degeneration from $\bbF_1$ to $Bl_p \bP(1,1,4)$ based on \cite[Section 5]{ADL19} and study the walls crossings on the surface pair $(Bl_p \bP(1,1,4),C)$. Let $[x_0,\cdots,x_3]$ be the homogeneous coordinate of $\bP(1,1,1,2)$ and then
consider the embedding \[ \bP^2_{[x,y,z]} \hookrightarrow \bP(1,1,1,2)\] 
given by 
\begin{equation}\label{P2em}
      x_0=x,\  x_1=y,\ x_2=z,\ x_3=xz-y^2
\end{equation}
and the the embedding \[ \bP(1,1,4)_{[x,y,z]} \hookrightarrow \bP(1,1,1,2)\] 
given by 
\begin{equation}\label{P114em}
   x_0=x^2,\ x_1=xy,\ x_2=y^2,\ x_3=z.
\end{equation}
It is easy to see that the hypersurface 
\begin{equation*}
    \cY:=\{ (s,[x_0,\cdots,x_3])\ |\ s \cdot x_3=x_0x_2-x_1^2\} \subset \bA^1 \times \bP(1,1,1,2)
\end{equation*}
under the natural projection $\cY \rightarrow \bA^1$ provides a  degeneration of $\bP^2$ to $\bP(1,1,4)$. Now we take a section $\gamma$ and then a blowup $\fX:=Bl_\gamma(\bA^1) \cY \rightarrow \cY$ along the image of $\gamma$. Then we get a family of  surface\begin{equation}\label{degeneration}
    \varphi: \fX \rightarrow \bA^1
\end{equation} such that $\fX_s \cong \bbF_1$ for $s \neq 0$ and $\fX_0 \cong Bl_p\bP(1,1,4)$. We denote the curve 
$C_0=2H_z+Q_0$  whose image is $B_0:=\pi(C_0)=\{z^3+z^2x^4=0\}$ where $Q_0$ is proper transform of curve $\{z+x^4=0\}$.  Let \[ B^-:=\{(xz-y^2)^2(xy-y^2+x^2)=0\}=2Q_1+Q_2 \] be the union of two plane conics $Q_1$ and $Q_2$ with multiplicity $2$ and $1$ and $C^-=:\pi^\ast B^--2E$.


\begin{pro}\label{2nddivc}
At the  wall $w=\frac{29}{106}$, there are  biratinal morphisms \[ \overline{P}^K_{\frac{29}{106}+\epsilon}\  \xrightarrow{p_+}\  \overline{P}^K_{\frac{29}{106}}  \  \xleftarrow{p_-}  \  \overline{P}^K_{\frac{29}{106}-\epsilon}\]  
where $p^-$ is an isomorphism and $p^+$ is a Kirwan type blowup. 
The exceptional locus $E_{w}^-$ is a point  parametrizing K-polystable pairs $(\bbF_1,(w-\epsilon)C^-)$ and $E_{w}^+$ is a divisor parametrizing K-polystable pairs $(Bl_p \bbF(1,1,4), C^+)$ described in first row of Table \ref{tabu}.  Moreover, $E_{w}^+$  is birational to the unigoal divisor $H_u$ in $\cF$.
\end{pro}
\begin{proof}
The $\bG_m$-action  with weight $(1,0,4)$  on $\bP(1,1,4)$ is induced from the $\bG_m$-action $\lambda$ on $\bP(1,1,1,2)$ \[ t\cdot [x_0,x_1,x_2,x_3]=[t^2x_0,tx_1,x_2,t^4x_3]\] under the embedding $\bP(1,1,4) \hookrightarrow \bP(1,1,1,2)$ via (\ref{P114em}). In particular $\lambda$ induces a 
 $\bG_m$-action on $\cY$ by \[ t\cdot (s,[x_0,x_1,x_2,x_3])=(t^{-2}s,[t^2x_0,tx_1,x_2,t^4x_3]).\]  By taking a section \[ \gamma: \bA^1 \rightarrow \cY , \ s\mapsto (s,[1,0,0,0])\]equivariant with respect to such $\bG_m$ action, then the degeneration construction $\fX \rightarrow \bA^1$ in (\ref{degeneration}) is also  $\bG_m$-equivariant. Then  the pair $(\bbF_1, C^-)$ admits a  special degeneration to the pair $(Bl_p \bbF(1,1,4), C_0)$, whose moduli point is  exactly the center $Z_w$ at the wall $w=\frac{29}{106}$ by part (1) of theorem \ref{mainthm1}. For the general curve $B^+$
of the form 
\begin{equation}\label{unigonalequ}
   z^3+z^2x^4+z^2yf_3(x,y)+zyf_7(x,y)+yf_{11}(x,y)=0 , 
\end{equation}
there is  a degeneration from $B^+$ to $B_0$ under the 1-PS $\lambda$ with weight $(1,0,4)$ in the Table \ref{Kwall2}. This will induce a 
$\bG_m$ degeneration from the pair $(Bl_p \bbF(1,1,4), C^+)$ to $(Bl_p \bbF(1,1,4), C_0)$ .  Thus  $(Bl_p \bbF(1,1,4), C^+)\in E^+_w$ by local VGIT interpretation of $\overline{P}^K_w$ .   Observe that such general curve $B^+$ in equation (\ref{unigonalequ}) is exactly general member in $|-2K_{Bl_p \bbF(1,1,4)}|$. By the dimension counting, \[\dim E^+=\dim |-2K_{Bl_p \bbF(1,1,4)}| \q \aut(Bl_p \bbF(1,1,4)) =\dim 17 .\] So $\dim E^-=0$ by dimension formula (\ref{dimvgit}). This proves $p^-$  must be an isomorphism and $p^+$ is a weighted blowup. By the computation of period point for K3 surfaces in \cite[Section 4.1.1]{PSW1} obtained by the pairs in  $E^+_w$,  $E^+_w$ is birational to $H_u$. Then we finish the proof. 
\end{proof}
\begin{rem}
    Under the transform $z\mapsto z-yf_3(x,y)$, the equation (\ref{unigonalequ}) has the normal form 
    \[ z^3+z^2x^4+zyf_7(x,y)+y^2f_{10}(x,y)=0.  \]
Then $\bG_m$ will act on the projective space $\bP^{18}\cong \bP(\bC[x,y]_7\oplus \bC[x,y]_{10})$  induced from the action on the  equation (\ref{unigonalequ}) and $E^+\cong \bP^{18}\q \bG_m$.  
\end{rem}

\subsubsection{The remaining walls on $Bl_p\bP(1,1,4)$ }
\begin{pro}\label{2ndflip}
At the third wall $w\in \{\frac{31}{110}, \frac{2}{7},\frac{35}{118} \}$, there are flips 
\[ \overline{P}^K_{w+\epsilon}\  \xrightarrow{p_+}\  \overline{P}^K_{w}  \  \xleftarrow{p_-}  \  \overline{P}^K_{w-\epsilon}\] 
The centers of $p_+$ and $p_-$ are given by Table \ref{Kwall2}. Moreover, 
\begin{enumerate} 
    \item  The exceptional locus  $E_w^{-}=p_-^{-1}(Z_w)$ parametrizes   $S$-equivalence classes of K-semistable pairs $(X,(w-\epsilon)\cdot C)$ where $B=\pi(C) \subset \bP(1,1,4)\subset \bP(1,1,1,2)$ or $B=\pi(C)\subset \bP^2 \subset \bP(1,1,1,2)$ is a  curve of complete intersection the form listed in the Table \ref{P-}. 
    \item  The exceptional locus  $E_w^{+}=p_+^{-1}(Z_w)$ parametrizes  $S$-equivalence classes of K-semistable  pairs $(X,(w+\epsilon)\cdot C)$ with  $B=\pi(C) \subset \bP(1,1,4)$ listed in the Table \ref{tabu}. 
    In addition,  $E_w^+$ is birational to a Noether-Lefschetz locus
\end{enumerate}
\begin{center}
\begin{table}[ht]
    \centering
      \begin{tabular}{ |c |c  |c|}
    \hline
      $c$ &   equation for  curve  $B$ on  $\bP(1,1,1,2)$  &   $\dim E^-$ \\ \hline 
      $ \frac{31}{110}-\epsilon$  & \makecell[c]{  $a_1x_3=x_0x_2+x_1^2$ \\  $x_3^3+x_3x_1x_0^3+a_2x_3^2x_0^2=0$ }   & $2$  \\ \hline 
      $ \frac{2}{7}-\epsilon$  & \makecell[c]{  $a_1x_3=x_0x_2+x_1^2$ \\  $x_3^3+x_2x_0^5+a_2 x_3x_1x_0^3+a_3x_3^2x_0^2 $}  & $3$  \\ \hline 
      $\frac{35}{118}-\epsilon$ &\makecell[c]{  $a_1x_3=x_0x_2+x_1^2$  \\  $x_3^3+x_3x_2x_0^3+x_1^2x_0^3+x_0^4f_2(x_1,x_2)=0$ }    & $4$ \\ \hline
   \end{tabular}
    \caption{Equations for  curves  in  exceptional locus $E^-$}
    \label{P-}
\end{table}
\end{center}
\end{pro}
\begin{proof}
We give a proof for $w=\frac{31}{110}$. The remaining cases are the same arguments and we leave the details  to the interested readers. 
Under the embedding (\ref{P114em}),  the curve $B \in Z_{w}$ in is cut out by the equations 
\begin{equation*}
    q_0:=x_0x_2-x_1^2,\  F_0=x_3^3+x_3x_1x_0^3=0.
\end{equation*}
Note that by Proposition \ref{degeneration}, t is shown for any K-semistable pairs $(X,C)$,  there is a blowup morphism $\pi: X \rightarrow \bP^2$ or $\pi: X \rightarrow \bP(1,1,4)$.  Let $\bC[x_0,x_1,x_2,x_3]_i$ be the space of homogeneous polynomial degree $i$ in  $x_0,x_1,x_2,x_3$, then the pair of polynomials $(q,f \mod q ) $ form a parameter space for  K-semistable pairs $(X,C)$
where  $q\in  \bC[x_0,x_1,x_2,x_3]_2$ and $f \in \bC[x_0,x_1,x_2,x_3]_6 $. Let $\lambda$ be the action on $\bP(1,1,1,2)$ by \[  t\cdot [x_0,x_1,x_2,x_3]=[t^4x_0,t^2x_1,x_2,t^7x_3] \]
induced from the 1-PS in the Table \ref{Kwall2}, which is the stabilizer group of the center $Z_w$.  Then the action $\lambda$ on the pair is 
\[ t\cdot (q,f\mod q):=(t^{-4} (t\cdot q), t^{-21} (t\cdot f))\]
where $t\cdot q$ or $t\cdot f$ is the action of $\lambda$ on $\bC[x_0,x_1,x_2,x_3]_2$ or $\bC[x_0,x_1,x_2,x_3]_6 $.  Then by the local VGIT interpretation of K-moduli at $\overline{P}^K_w$,  $E^-$ also parametrizes pairs $\{ (q, f \mod q )\}$ such that the curve defined by $\{q=0,\ f=0\}$ is nodal at $[1,0,0,0]$ at least and  \[ \mathop {\lim} \limits_{t \rightarrow 0}  t(q,f \mod q):=(\mathop {\lim} \limits_{t \rightarrow 0}  t^{-4} \cdot tq, \mathop {\lim} \limits_{t \rightarrow 0}  t^{-21} \cdot tf\mod q )=(q_0,F_0\mod q) \] In this way, we obtain equation 
\begin{equation*}
    q=q_0+x_0l(x_0,x_1)+a_1x_3, \ \ f=F_0+a_2x_3^2x_0^2
\end{equation*}
 for curves $B\in E_w^-$. Up to the possible action under $x_3\mapsto x_3+f_2(x_0,x_1,x_2)$, the normalised equation is obtained for $B$ as in Table \ref{P-}. This proves the description in (1).  By the parallel computation and discussion for curves in $E_w^+$, we prove the the  description in (2). 
\end{proof}


\subsection{Proof of part (2) of theorem \ref{mainthm1}}
This is just the combination of Proposition \ref{1stdivc} ,  Proposition \ref{1stflip} and  Proposition \ref{2nddivc},  Proposition \ref{2ndflip}. 

\section{K-moduli v.s. HKL } \label{sect6}
In in section, we established the relation of K-moduli spaces $\overline{P}^K_c$ with the Hassett-Keel-Looignega(HKL) program for the moduli space of lattice polarised K3 surfaces studied in \cite{PSW1}. 
\begin{lem}\label{normaility}
The K-moduli space $\overline{P}^K_c$ is a normal projective variety for any $\frac{1}{14}<c<\frac{1}{2}$. 
\end{lem}
\begin{proof}
By the structure of K-moduli and Luna's slice theorem, it is known that for any pairs $[X,C]\in \overline{P}^K_c$ (see \cite[Theorem 3.33]{ADL19}), there is a \'etable map 
\begin{equation*}
    U \q \aut(X,C)\  \rightarrow \overline{P}^K_c
\end{equation*}
whose image is a open neighborhood of  $[X,C]$. As $\aut(X,C)$ is a reductive group, then it is sufficient to show the deformation of $\bQ$-Gorenstein del Pezzo pair $(X,C)$ has no obstruction so that $U$ can be choose as a open subset of the 1st order deformation space of the  pair $(X,C)$, which is smooth and thus the GIT $U \q \aut(X,C)$ is normal by the general results of Mumford's GIT \cite[Charter 0 \S 02]{MKGIT}. 

Now let us compute the 1st order deformation space and obstruction space to finish the proof.  Recall from  \cite[Chapter 3]{Sernesi06},  deformation space and obstruction space of a klt pair $(X,C)$ in dimension $2$ are given by 
\[\mathrm{Def}_{(X,C)}=\Ext^1(\Omega^1_X(\log (C)_{red}),\cO_X),\ \  \mathrm{Obs}_{(X,C)}=\Ext^2(\Omega^1_X(\log  (C)_{red}),\cO_X). \]
where $\Omega^1_X(\log (C)_{red})$ is the logarithmetic differential sheaf along the reduced part $(C)_{red}$ of the curve $C$. 
By \cite[Proposition 3.1]{HP10}, \[H^2(X,T_X)\cong \Ext^2(\Omega^1_X(\log C),\cO_X)=0\]  by Serre duality.  Note that there is the short exact sequence 
\begin{equation}\label{sess}
    0 \rightarrow  \Omega^1_X  \rightarrow  \Omega^1_X(\log C) \rightarrow \mathop{\oplus} \limits_{i} \cO_{C_i} \rightarrow  0
\end{equation}
where $C_i$ is the irreducible component of $(C)_{red}$. Thus by taking the $\Ext(-,\cO_X)$ of (\ref{sess}), it is sufficient to show the vanishing $\Ext^2( \cO_{C_i},\cO_X)=0$ for each component $C_i$.  By taking  $\Ext(-,\cO_X)$ for the short exact sequence 
\[ 0 \rightarrow \cO_X(-C_i) \rightarrow \cO_X \rightarrow \cO_{C_i} \rightarrow 0 ,\]
it is easy to see that the vanishing $\Ext^2( \cO_{C_i},\cO_X)=0$ will follows from the following  vanishings 
\begin{equation}\label{vanishing}
\begin{split}
   & \Ext^2( \cO_X,\cO_X)=H^2(X,\cO_X)=0,\\
   & \Ext^1( \cO_X(-C_i),\cO_X)=H^1(X,\cO_X(C_i))=0.
\end{split} 
\end{equation}
 The first vanishing result in  (\ref{vanishing}) is just Kodaira vanishing for del Pezzo surface. The second  in  (\ref{vanishing}) will follow from the proof of \cite[Lemma2.1]{MS20} since $C_i$ is an effective divisor on the klt del Pezzo surface.  
\end{proof}
Recall there is a uniform embedding into $\bP^N$ for all K-semistable pair $(X,D)$ and denote $\Hilb_i=\Hilb_i(\bP^N)$ the Hilbert scheme  with Hilbert polynomial 
 \[\chi_1(t)=\frac{m^2d}{2} t^2+\frac{md}{2} t+1, \  \ \  \chi_2(t)=2md t-d.\]
Let $\pi_c: (\fX,\cD;\cL) \rightarrow  Z_c^{red} $ be the universal family of polarised del Pezzo pairs of degree $8$ where $ Z_c^{red}  \subset \Hilb_1 \times  \Hilb_2$ is the reduced locally closed subscheme in the product of Hilber scheme  parametrizing the incidence pair. By  \cite{ADL19}, the K-moduli $\overline{P}^K_c$ is the good moduli space of quotient stack $[Z_c^{red}/\PGL(N+1)]$. Denote $\lambda_{c,c'}$ the descent of the log CM line bundle with coefficient $c'$ on $\overline{P}^K_c$, that is, $\lambda_{c,c'}$ is the descent of 
\[ \pi_{c \ast} (-K_{\pi_c}-c' \cD)^3.\] 
It is known that $\lambda_c=\lambda_{c,c}$ is an ample $\bQ$-line bundle on $\overline{P}^K_c$ 
 for all $c\in (\frac{1}{14},\frac{1}{2})\cap \bQ$ by theorem \ref{cmample}. 
Moreover, there is birational contraction map \[ \overline{P}^K_{\frac{1}{2}-\epsilon} \dashrightarrow \overline{P}^K_c \] 
for any $c $. Then
arguments as in \cite[Theorem 9.4]{ADL19} will show 
\begin{equation}\label{kinterporation}
     \overline{P}^K_c\cong \Proj \big( R( \overline{P}^K_{\frac{1}{2}-\epsilon}, \lambda_{\frac{1}{2}-\epsilon,c}) \big)
\end{equation}
for $c\in (\frac{1}{14},\frac{1}{2})\cap \bQ$  and $0 <\epsilon \ll 1$.

\begin{thm}\label{thm-identification}
There is a natural isomorphism of projective varieties
$\overline{P}^K_c \cong \cF(s)\ $ induced by the period map   under the transformation  \[ \ s=s(c)= \frac{1-2c}{56c-4}\] where $\frac{1}{14}<c<\frac{1}{2}$.
\end{thm}

\begin{proof}
Let $p:\overline{P}^K_{\frac{1}{2}-\epsilon} \dashrightarrow \cF^\ast$ be the period map. By  the results in Section 5, we know 
\[  p^{-1} : \cF \hookrightarrow  \overline{P}^K_{\frac{1}{2}-\epsilon} \]
is a open immersion whose image is a big open subset of  $\overline{P}^K_{\frac{1}{2}-\epsilon}$. 
Thus by the normality of K-moduli space $  \overline{P}^K_c$ proved in Lemma \ref{normaility} and (\ref{kinterporation}), it is enough to show the pullback $(p^{-1})^\ast \lambda_{\frac{1}{2}-\epsilon,c}$ on $\cF^\ast$ is proportional  to
\[ \lambda + \frac{1-2c}{56c-4}( H_h+ 25 H_u) .\]
By interpolation formula of CM line bundles in \cite[Proposition3.35]{ADL19}, we have
\begin{equation} \label{unmcm}
    (1-2c)^{-2}\cdot \lambda_{\frac{1}{2}-\epsilon,c}=(1-2c) \cdot \lambda_{\frac{1}{2}-\epsilon,0}+48c \cdot \lambda_{\frac{1}{2}-\epsilon,Hdg}.
\end{equation}
Then we claim 
\begin{equation}\label{claimcm}
    p^{-1\ \ast}\lambda_{\frac{1}{2}-\epsilon,Hdg}= \lambda,\ \  p^{-1\ \ast} \lambda_{\frac{1}{2}-\epsilon,0}= H_h+25H_u-4\lambda
\end{equation}
Assume the claim (\ref{claimcm}), the formula (\ref{unmcm}) will imply 
\begin{equation*}
   \begin{split}
        R( \overline{P}^K_{\frac{1}{2}-\epsilon}, \lambda_{\frac{1}{2}-\epsilon,c} \big)=R(\cF^\ast, (56c-4)\lambda+(1-2c)(H_h+25H_u)).
   \end{split} 
\end{equation*} 
Thus, to prove the theorem it is sufficient to prove the claim (\ref{claimcm}). The first identity is obtained by adjunction as in the proof of \cite[Theorem 6.2]{ADL22}.  As $\pic(\cF)_\bQ$ is generated by $\lambda, H_h,H_u$, we may assume 
\[ p^{-1\ \ast} \lambda_{\frac{1}{2}-\epsilon,0}=a_h H_h+a_uH_u+a_\lambda \lambda,\ a_u,a_h,a_\lambda\in \bQ . \]
By Proposition \ref{1}, the section ring \[R(c):=R(\cF, 48c \lambda+(1-2c) \cdot (a_h H_h+a_uH_u+a_\lambda \lambda)\] will satisfy 
\begin{equation*}
  R(c)= \begin{cases}
       0 ,\ & \ 0<c<\frac{1}{14}\\
       \bC, \ & \ c=\frac{1}{14}.
  \end{cases} 
\end{equation*}
 Since $\lambda$ is ample while $H_h$  and $H_u$ are contractable,  the coefficient of $\lambda$ in $48c \lambda+(1-2c) \cdot (a_h H_h+a_uH_u+a_\lambda \lambda)$ must vanish, that is, \[\frac{48}{14}+\frac{12}{14}a_\lambda =0 . \] This shows $a_\lambda=-4$. In particular, the CM line bundle $\lambda_{\frac{1}{2}-\epsilon,c}$ is proportional to
 \[ \lambda+\frac{1-2c}{56c-4}\cdot (a_h H_h+a_uH_u). \]
By the computation in \cite[Theorem 1.2]{PSW1}, 
\begin{equation}\label{restrction}
    (\lambda+H_u)|_{H_u}=0,\  \ (\lambda+H_h)|_{H_h}=0.
\end{equation} 
where we still use $H_u$ and $H_h$  denote their birational transform in the K-moduli space $\overline{P}^K_c$. Note that by the computation in Proposition \ref{1stdivc} and Proposition \ref{2nddivc}, the walls where the proper transforms of $H_h$ and $H_u$ in $P_c^K$ appear as divisorial contractions  are $w_h=\frac{5}{58}$ and $w_u=\frac{29}{106}$ . Then  we have
 \begin{equation*}
     \begin{split}
       \frac{1-2w_h}{56 \cdot w_h-4} \cdot  a_h=1,  \ \   \frac{1-2w_u}{56 \cdot w_u-4} \cdot  a_u=1
     \end{split}
 \end{equation*}
 as in the proof of  \cite[Theorem 6.2]{ADL22}. It turns out that $a_h=1,\ a_u=25$. Then we finish the proof.
\end{proof}
\begin{rem}
    Under the the transformation $s=s(c)= \frac{1-2c}{56c-4}$, the walls   in  (\ref{mainthm1}) on hyperelliptic divisor $H_h$ will be 
    \begin{equation*}
        \begin{split}
            s\ \in \ &W_{h,A}=\{1,  \frac{1}{2},  \frac{1}{3},   \} \cup W_{A'}= \{ \frac{1}{4},\frac{1}{6},\frac{1}{8}\}  \\ 
            & \cup W_D=\{\frac{1}{4},\frac{1}{6},\frac{1}{8}, \} \cup W_{D'}=\{ \frac{1}{10},\frac{1}{12} , \frac{1}{16} \} \cup W_E=\{\frac{1}{10},\frac{1}{16},\frac{1}{28} \}
        \end{split}
    \end{equation*}
   and the walls on unigonal divisor $H_u$ will be 
    \[s\in W_u=\{\frac{1}{25},\frac{1}{27},\frac{1}{28},\frac{1}{31}\}.\]
    These numbers are just walls from prediction of arithmetic side in \cite[Prediction 1 and Prediction 2]{PSW1}. In particular, we verify the walls in conjecture \ref{HKL8}. 
\end{rem}

\section{Further discussion}\label{sect7}
\subsection{Application to K-stability of Fano 3-fold pairs }
Let us recall the construction of log  Fano $3$-fold pair associated to the del Pezzo surface via the cone construction. Let $X$ be a normal surface such that $-K_X$ is ample and Cartier.  Denote $C \in |-2K_X|$ an effective curve. The $3$-fold $V$ is defined as projective cone
\[ V:=C_p(X,-K_X)=\Proj(\mathop{\oplus}\limits_{m\ge 0}\mathop{\oplus}\limits_{n=0}^m H^0(X,-nK_X))t^{m-n})\]
and the surface $S$ as an anti-canonical divisor on $V$  is obtained by  double covering $S \rightarrow X$  branched along the curve $C$. We call such  $(V,S)$   the log Fano 3-fold pair associated to the del Pezzo surface pair $(X,C)$.

A powerful results proved in \cite[Theorem 5.2]{ADL22} is the following 
\begin{thm}\label{kcone}
  Notation as above, then  the log Fano 3-fold pair $(V,\frac{4c+1}{3}S)$ is K-semistable if and only if the del Pezzo surface pair $(X, cC)$ is K-semistable.
\end{thm}
\begin{rem}
    The $3$-fold $V$ has volume $(-K_V)^3=8(-K_X)^2$. 
\end{rem}
As a application of the results in Section 6, we have  
\begin{cor} \label{3foldpair}
Let $(V,S)$  be the log Fano 3-fold pair obtained by the cone construction of del Pezzo pairs $(\bbF_1,C)$ for curve $C$ given by the Table \ref{Kwall1}, then 
 \[ kst(V,S)=\{\ c=\frac{11+n}{27+n} \ |\  n=1,2,\cdots, 5\ \} \cup \{\ c=\frac{3+n}{11+n}\ |\ n=6, 7,8,9,11\  \}.  \]
\end{cor}
\begin{proof}
    By the computation of walls in Section 6, the assertion is a direct consequence of theorem  \ref{kcone}.
\end{proof}

\subsection{KSBA moduli vs toroidal compactification}

\begin{defn}
    A KSBA stable pair $(X,\Delta)$ consists of a projective variety $X$ with an effective $\bQ$-divisor $\Delta$ such that 
    \begin{enumerate}
        \item $K_X+\Delta$ is ample,
        \item  $(X,\Delta)$  has slc singularities. That is, the pair $(X^\nu,\Delta^\nu)$ is lc  where $\nu: X^\nu \rightarrow X$ is a normalization of $X$ such that 
        \[ K_{X^\nu}+\Delta^\nu=\nu^\ast (K_X+\Delta). \]
    \end{enumerate}
\end{defn}
If we add coefficient $c>\frac{1}{2}$ to the boundary curve $C$ of  del Pezzo pair $(X,C)$, then $(X,cC)$ will be a KSBA stable pair. Due to the frame work of Kollar-Shepherd-Barron (see \cite{KSB88} or  Kollar's book \cite{kollar} for the details), there is an irreducible component $\overline{P}_{d,c}^{KSBA}$ of a complete projective scheme  which parametrizes  KSBA stable pairs such that the general member is the isomorphic class of $(X,cC)$.  Moreover, varying $c>\frac{1}{2}$ there is also wall-crossing phenomenon for KSBA moduli space by \cite[Theorem 1.1]{ABIP21}.

Given a stable pair  $(X,C)\in \overline{P}_{d,\frac{1}{2}+\epsilon}^{KSBA}$ and by the double cover construction $\phi: Y \rightarrow X$ branched along  curve $C$, one can obtain another  pair  $(Y,R)$ where  $K_Y=0$  (a singular K3 ) and $R=\phi^{-1}(C)$ is a $\bQ$-divisor. It turns out $(Y,c'R)$ is another KSBA stable pair  for any $c'>0$ by \cite[Proposition 4.1 ]{AEH}.  Under such construction Alexeev-Engel-Han 
 show (\cite[Theorem 4.2 and Corollary 4.3]{AEH})  $\overline{P}_{d,\frac{1}{2}+\epsilon}^{KSBA}$ is also the coarse moduli space of KSBA stable $(Y,\epsilon R)$ where $0< \epsilon \ll 1$. By realising $R$ as a recognizable divisor, they also show there is a normalization map 
\begin{equation*}
   \mu: \cF^{\Sigma} \rightarrow \overline{P}_{d,\frac{1}{2}+\epsilon}^{KSBA}
\end{equation*}
where $\cF^{\Sigma}$ is a certain semi-toric compactification of $\cF$ determined by a $O(\Lambda)$-admissble semifan $\Sigma$. From arithmetic compactification side, there is a natural contraction morphism from the  semi-toric compactification $\cF^{\Sigma}$ to the Baily-Borel compactification $\cF^{\ast}$
\begin{equation*}
    \cF^{\Sigma} \rightarrow  \cF^{\ast}=P_d^\ast.
\end{equation*}
 It maps the simple normal crossing boundary divisors $\cF^{\Sigma}-\cF=D_1\cup \cdots\cup  \cdot \cup D_m$ to $1$-dimensional cusps $B_1,\cdots,B_m$ of the boundaries $\cF^\ast-\cF$.
Let $\overline{P}_{d,\frac{1}{2}+\epsilon}^{KSBA,\circ}$ be the Zariski open subset parametrizing ADE stable pairs. We expect that the period map $p:\overline{P}_{d,\frac{1}{2}+\epsilon}^{KSBA,\circ} \rightarrow  \cF^\ast$  will extend to a morphism $p:\overline{P}_{d,\frac{1}{2}+\epsilon}^{KSBA} \rightarrow  \cF^\ast$ such that there is a commutative diagram
\begin{equation*}
    \begin{tikzcd}
      \cF^{\Sigma}  \arrow[dr] \arrow[rr,"\mu"]  &  &  \overline{P}_{d,\frac{1}{2}+\epsilon}^{KSBA} \arrow[dl,"p"]  \\
       & \cF^\ast &
    \end{tikzcd}
\end{equation*}
Then as analogue to HKL,  we consider the projective scheme
\begin{equation*}
   \cF^{\Sigma}(b_1,\cdots, b_m):=  \Proj \big (R(\cF^{\Sigma},\lambda+\mathop{\sum} \limits_{i=i}^m b_iD_i)\ \big), \ b_i \in [0,1) \cap \bQ 
\end{equation*}
As the pair $(\cF^{\Sigma},\mathop{\sum} \limits_{i=i}^m b_iD_i)$ is klt, then the general result of \cite[Corollary 1.1.2 ]{BCHM} will imply $\cF^{\Sigma}(b_1,\cdots, b_m)$ is a projective variety. An interesting question is 
\begin{question}
    If we vary the parameter $(b_1,\cdots, b_m)$, is it true that  $\cF^{\Sigma}(b_1,\cdots, b_m)$ is the normalization of KSBA moduli space $\overline{P}_{d,c}^{KSBA}$ for certain $c=c(b_1,\cdots, b_m)>\frac{1}{2}$ ?
\end{question}
If the answer is positive, then the similar arithmetic strategy in \cite[Section 10]{FLLST} should provide an effective algorithm to  find walls for KSBA moduli spcace. 
\subsection{ Wall crossing relating K-moduli and KSBA moduli}
Motivated by the  \cite[Conjecture 9.19]{ADL19},  we ask the following 
\begin{question}\label{Qlogcy}
  Is there a good stability notion for log CY surfaces $(X,\frac{1}{2}C)$  such that  the  log CY  surfaces $(X,\frac{1}{2}C)$ with such stability forms a good  moduli problem , which admits a good moduli space $\overline{P}_{d}^{CY}$ and  $\overline{P}_{d}^{CY}$ fits into the  wall-crossing of K-moduli $\overline{P}_{d,\frac{1}{2}-\epsilon}^{K}$ and KSBA-moduli $\overline{P}_{d,\frac{1}{2}+\epsilon}^{KSBA}$ ?
\end{question}

An evidence for a positive answer to the Question \ref{Qlogcy} is provided by the following result. 
\begin{pro}
    $\cF^\ast$ is the ample model of of the Hodge line bundle $\overline{P}^K_{\frac{1}{2}-\epsilon}$. 
\end{pro}
\begin{proof}
     The proof  is parallel  to that of  \cite[Theorem 6.5]{ADL19}.  Let $U \subset \cF^\ast$ be the open subset parametrizing  K3 surface with involution $(X,\tau)$ such that $X/\tau$ is either  $\bbF_1$ or  $\widetilde{Bl_p \bP(1,1,4)}$ the minimal resolution of $Bl_p \bP(1,1,4)$.  Since $\rho(\cF^\ast)=3$, then $U\subset \cF^\ast$ is a big open subset, i.e., $codim (\cF^\ast -U) \ge 2$. By the explicit wall-crossing description in  Section \ref{explicitwc}, $U \subset \overline{P}^K_{\frac{1}{2}-\epsilon}$ is also a big open subset. Denote $\lambda_{Hodge}$ the Hodge line bundle on $\overline{P}^K_{\frac{1}{2}-\epsilon}$.  Note that  $\lambda_{Hodge}|_U$ is obtained by the restriction of Hodge line bundle $\lambda$ on locally symmetric variety $\cF=\cD \setminus \Gamma$.  As arguments in the proof of  \cite[Theorem 6.5]{ADL19},  the line bundle $\lambda_{Hodge}$ is big and semiample on $\overline{P}^K_{\frac{1}{2}-\epsilon}$. Thus, 
     \begin{equation*}
         \begin{split}
             \cF^\ast &  \cong  \Proj(R(\cF, \lambda )) \cong  \Proj(R(U, \lambda|_U )) \\
 & \cong  \Proj(R(U, \lambda_{Hodge}|_U )) \cong  \Proj(R(\overline{P}^K_{\frac{1}{2}-\epsilon}, \lambda_{Hodge} )).
         \end{split}
     \end{equation*}
     This finishes proof.
\end{proof}
Recently, the authors in \cite{ABBDILW23} succeed in constructing the  moduli of log Calabi-Yau pairs in  plane curves case, which  connects  the KSBA moduli and the K-moduli for the plane curve pairs. Their work provides further evidence to the Question \ref{Qlogcy}. We expect extending their work to give an answer to this question. 
\begin{rem}
    The full explicit wall crossing description  of K-moduli $\overline{P}^K_{d,c}$ for other degree $d$ will imply the same result that $\cF^\ast$ is the ample model of of the Hodge line bundle $\overline{P}^K_{d, \frac{1}{2}-\epsilon}$. 
\end{rem}

\vspace{0.5cm}
\bibliographystyle{alpha}
\bibliography{references}
\end{document}